\documentclass{amsart}
\usepackage{amssymb}
\usepackage[v2,cmtip]{xy}

\chardef\bs=`\\ 
\newcommand{\myop}[1]{\mathop{\textstyle #1}\limits}
\newcommand{\ssdot}{\bullet}%{{\scriptscriptstyle\bullet}}
\newcommand{\subdot}{_\ssdot}
\newcommand{\supdot}{^\ssdot}
\newcommand{\op}{\mathrm{op}}
\newcommand{\num}[1]{\##1}

\newcommand{\Q}{\bQ}
\newcommand{\Z}{\bZ}
\newcommand{\ZP}{\Z/p\Z}
\newcommand{\Zhat}{\Z^{\wedge}}
\newcommand{\ZPhat}{\Zhat_{p}}
\newcommand{\Qhat}{\Q^{\wedge}}
\newcommand{\FP}{\bF_{p}}
\newcommand{\FPbar}{{\bar\bF}_{p}}
\newcommand{\W}{\bW}
\newcommand{\card}{\varkappa}
\newcommand{\continuum}{\mathfrak{c}}
\newcommand{\Zprime}{Z'}

\newcommand{\loopf}{\Lambda_{f}}
\newcommand{\rat}[1]{#1_{\Q}}
\newcommand{\finite}[1]{#1^{\wedge}}
\newcommand{\Wmap}[1]{#1_{\scriptscriptstyle \W}}
\newcommand{\Qmap}[1]{#1_{\scriptscriptstyle \Q}}

\newcommand{\Ei}{\mathfrak{E}}
\newcommand{\SP}{\mathfrak{Top}}
\newcommand{\SSet}{\oS}
\newcommand{\Ecard}{\Ei'}
\newcommand{\Scard}{\SSet'}
\newcommand{\Modcat}{\mathfrak{M}}

\def\HoSSMap(#1,#2){[#1,#2]}
\def\HoSPMap(#1,#2){[#1,#2]}
\def\HoEiMap(#1,#2){[#1,#2]_{\Ei}}

\def\HoSP{Ho(\SP)}

\newcommand{\retr}{\epsilon}
\newcommand{\adj}{U}
\newcommand{\deradj}{\mathbf{U}}
\newcommand{\Tder}{\mathbf{T}}
\newcommand{\AQS}{D}

\newcommand{\rfib}{\rto|>>\tip}
\newcommand{\lfib}{\lto|>>\tip}
\newcommand{\dfib}{\dto|>>\tip}
\newcommand{\rcof}{\rto|<<\tip}
\newcommand{\lcof}{\lto|<<\tip}

\newtheorem{thm}{Theorem}[section]
\newtheorem{cor}[thm]{Corollary}
\newtheorem{lem}[thm]{Lemma}
\newtheorem{prop}[thm]{Proposition}
\newtheorem*{main}{Main Theorem}
\newtheorem{theorem}{Theorem}

\theoremstyle{definition}
\newtheorem{defn}[thm]{Definition}
\newtheorem{notn}[thm]{Notation}

\theoremstyle{remark}

\numberwithin{equation}{section}
\newcommand{\term}[1]{\textsl{#1}}
\makeatletter
\let\c@equation\c@thm
\makeatother

\newcommand{\bE}{{\mathbb{E}}}
\newcommand{\bF}{{\mathbb{F}}}
\newcommand{\bQ}{{\mathbb{Q}}}
\newcommand{\bW}{{\mathbb{W}}}
\newcommand{\bZ}{{\mathbb{Z}}}

\let\opsymbfont\mathcal 
\newcommand{\oD}{{\opsymbfont{D}}}
\newcommand{\oE}{{\opsymbfont{E}}}
\newcommand{\oS}{{\opsymbfont{S}}}

\newcommand{\iso}{\cong}     % preferred isomorphism symbol
\newcommand{\homeq}{\simeq}  % homotopy symbol

\newcommand{\einf}{\text{$E_\infty$} } % should never be at end of sentence
\renewcommand{\to}{\longrightarrow}
\newcommand{\from}{\longleftarrow}

\newcommand{\overto}[1]{\xrightarrow{\,#1\,}}

\def\quickop#1{\expandafter\DeclareMathOperator\csname #1\endcsname{#1}}
\quickop{id} \quickop{Lim} \quickop{Colim} \quickop{Gal}\quickop{diag}

\begin{document}

\title{Cochains and Homotopy Type}

\author{Michael A. Mandell}
\address{Department of Mathematics, University of Chicago, Chicago, IL}
\email{mandell@math.uchicago.edu}
\thanks{The author was supported in part by NSF postdoctoral
fellowship DMS 9804421.}
\date{\today}

\begin{abstract}
Finite type nilpotent spaces are weakly equivalent if
and only if their singular cochains are quasi-isomorphic as \einf
algebras.  The cochain functor from the homotopy category of finite
type nilpotent spaces to the homotopy category of \einf algebras is
faithful but not full.

\end{abstract}

\maketitle

\section*{Introduction}

Motivating questions in algebraic topology often take the form of
finding algebraic invariants for some class of topological spaces that
classify those spaces up to some useful equivalence relation, like
homeomorphism, homotopy equivalence, or even something weaker.  The
fundamental work of Quillen \cite{quilrat} and Sullivan \cite{sullrat}
describes a form of algebraic data that classifies simply connected
spaces up to rational equivalence, the equivalence relation generated
by maps of spaces that are rational homology isomorphisms.  In
particular, the latter work associates to each space a rational
commutative differential graded algebra (CDGA) called the polynomial
De\,Rham complex, which is closely related to the algebra of
differential forms when the space is a smooth manifold and closely related to
the singular cochain complex in general.  Simply connected spaces of finite type
(homology finitely generated in each degree) are rationally equivalent
if and only if their associated polynomial De\,Rham complexes are
quasi-isomorphic CDGAs.

More recently, in \cite{einf} the author proved an analogous theorem for
$p$-equivalence, the equivalence relation generated by maps that are
$\ZP$ homology isomorphisms.  Here the relevant sort of algebra is an
\einf algebra, an up-to-homotopy generalization of a commutative
differential graded algebra.  The singular cochain complex of a space
with coefficients in a commutative ring $R$ has the natural structure
of an \einf $R$-algebra.  Simply connected spaces of finite type are
$p$-equivalent if and only if their cochains with coefficients in
$\FP$ ($=\ZP$) are quasi-isomorphic \einf $\FP$-algebras.  As
discussed below, coefficients in the algebraic closure $\FPbar$ lead
to even stronger results.

The purpose of this paper is to prove that the singular cochains with
integer coefficients, viewed as \einf algebras, classify finite type
simply connected spaces up to weak equivalence.  We prove the
following theorem.

\begin{main}
Finite type nilpotent spaces $X$ and $Y$ are weakly equivalent if and
only if the \einf algebras $C^{*}(X)$ and $C^{*}(Y)$ are quasi-isomorphic.
\end{main}

Weak equivalence is the equivalence relation generated by the maps of
spaces that induce bijections of components and isomorphisms of all
homotopy groups for all base points; such maps are themselves called weak
equivalences.  By a theorem of Whitehead, spaces homotopy equivalent
to CW complexes (which include most spaces of geometric interest) are
weakly equivalent if and only if they are homotopy equivalent.
Nilpotent is a generalization of simply connected: It means that the
fundamental group at each point is nilpotent and acts nilpotently
on the higher homotopy groups.

The invariants discussed in \cite{einf,quilrat,sullrat} do more than
distinguish equivalence classes of spaces, they distinguish
equivalence classes of maps.  We can consider the category obtained
from the category of spaces by formally inverting the rational
equivalences; this category is called the rational homotopy category.
Likewise, we can form the homotopy category of CDGAs by formally
inverting the quasi-isomorphisms.  Then \cite{sullrat} proves that the
full subcategory of the rational homotopy category of finite type
nilpotent spaces embeds as a full subcategory of the homotopy category
of CDGAs.  Similarly \cite{einf} proves that the full subcategory of
the $p$-adic homotopy category ($p$-equivalences formally inverted) of
finite type nilpotent spaces embeds as a full subcategory of the homotopy
category of \einf $\FPbar$-algebras. For the homotopy category of
$\FP$-algebras, the functor is faithful but not full.

In the integral situation, the category obtained from spaces by
formally inverting the weak equivalences is called simply the homotopy
category; it is equivalent to the category of CW complexes and
homotopy classes of maps.  We show that the singular cochain functor
with integer coefficients, viewed as a functor from the homotopy
category to the homotopy category of \einf algebras is faithful but is
not full, just as in the case for coefficients in $\FP$.  First, for
faithfulness, we prove the following refinement of the Main Theorem.
In it, $\HoSPMap(X,Y)$ denotes the set of maps from the space $X$ to
the space $Y$ in the homotopy category, and $\HoEiMap(A,B)$ denotes
the set of maps from the \einf algebra $A$ to the \einf algebra $B$ in
the homotopy category of \einf algebras.

\begin{theorem}
There is a function $\retr\colon \HoEiMap(C^{*}Y,C^{*}X)\to
\HoSPMap(X,Y)$, natural in spaces $X$ and finite type nilpotent spaces
$Y$, such that the composite
\[
\HoSPMap(X,Y)\overto{C^{*}}
\HoEiMap(C^{*}Y,C^{*}X)\overto{\retr}
\HoSPMap(X,Y)
\]
is the identity.  If in addition $X$ is finite type, then for any $f$
in $\HoEiMap(C^{*}Y,C^{*}X)$, the maps 
$f$ and $C^{*}(\retr(f))$ induce the same map on cohomology.
\end{theorem}

We suspect that last the statement is true for arbitrary $X$, but our
techniques do not appear to suffice.

The significance of Theorem~A is that even though $C^{*}$ is not a full
embedding, it is still sufficient for many classification problems about
existence and uniqueness of maps.  For example, for any finite type
space $X$ and any finite type nilpotent space $Y$, Theorem~A implies:
\begin{enumerate}
\item (Existence) Given any $\phi$ in $\HoEiMap(C^{*}Y,C^{*}X)$, there
exists $f$ in $\HoSPMap(X,Y)$ such that $H^{*}f=H^{*}\phi$. 
\item (Uniqueness) For $f,g$ in $\HoSPMap(X,Y)$, $f=g$ if and only if
$C^{*}f=C^{*}g$ in the homotopy category of \einf algebras.
\end{enumerate}
We should also note that the Main Theorem is an immediate consequence
of Theorem~A: By Dror's generalization of the Whitehead Theorem
\cite{drorwhitehead}, a map between nilpotent spaces is a weak
equivalence if and only if it induces an isomorphism on homology.  For
finite type spaces, a map induces an isomorphism on homology if and
only if it induces an isomorphism on cohomology.

Finally, to show that the functor $C^{*}$ is not full, we describe the set of
maps $\HoEiMap(C^{*}Y,C^{*}X)$ for $Y$ finite type nilpotent.  As we
explain in the next section, the functor $C^{*}$ has a contravariant
right adjoint $\deradj$.  That is, there exists a functor $\deradj$
from the homotopy category of \einf algebras to the homotopy category
and a natural bijection
\[
\HoSPMap(X,\deradj A)\iso \HoEiMap(A,C^{*}X)
\]
for any space $X$ and any \einf algebra $A$.  Thus, to understand
$\HoEiMap(C^{*}Y,C^{*}X)$, we just need to understand the space
$\deradj C^{*}Y$.  Let $\finite{Y}$ denote the finite completion of
$Y$ and let $\rat{Y}$ denote the rationalization of $Y$, as defined
for example in \cite{sullght}.  Let $\Lambda\finite{Y}$ denote the
free loop space of $\finite{Y}$, and consider the map $\Lambda \finite{Y}\to
\finite{Y}$ given by evaluation
at a point.  Define $\loopf Y$ to be the
homotopy pullback of the map $\Lambda \finite{Y}\to
\finite{Y}$ along the finite completion map $Y\to \finite{Y}$.
\[
\diagram
\loopf Y\rdotted|>\tip\ddotted|>\tip&\Lambda \finite{Y}\dto\\
Y\rto&\finite{Y}
\enddiagram
\]
We prove the following theorem.

\begin{theorem}
Let $Y$ be a finite type nilpotent space.  Then $\deradj C^{*}Y$ is weakly
equivalent to $\loopf Y$.
\end{theorem}

We have that $\pi_{n}\loopf Y=\pi_{n}Y \oplus \pi_{n+1}\finite Y$ at
each base point of $\loopf Y$.  When $Y$ is finite type nilpotent
$\pi_{n+1}\finite{Y}$ is the pro-finite completion of $\pi_{n+1}Y$.
When $Y$ is connected and simply connected, $\loopf Y$ is connected
and simple (the fundamental group is abelian and acts trivially on the
higher homotopy groups), and so $\HoSPMap(S^{n},\loopf Y)\iso
\pi_{n}\loopf Y$.

It follows from Theorem~B that $C^{*}\colon \HoSPMap(X,Y)\to
\HoEiMap(C^{*}Y,C^{*}X)$ is usually not surjective.  As a concrete
example, when we take $X=S^{2}$ and $Y=S^{3}$, the set $\HoSPMap(X,Y)$
is one point and the set $\HoEiMap(C^{*}Y,C^{*}X)$ is uncountable.

\bigskip

It seems plausible that the finite type hypotheses may be dropped by
considering a category of (flat) \einf coalgebras in place of \einf
algebras.  Although some authors have claimed this sort
of result \cite{smirnov} or better results \cite{smithcoalg}, as of
this writing no one has found a correct proof.

\subsection*{\bf Acknowledgments}
It is a pleasure for the author to thank a number of people for their
interest in and contribution to this work and \cite{einf}, where the
acknowledgments were inadvertently omitted.  This paper could not have
been written without the aid of many useful conversations with Mike
Hopkins or without the advice and encouragement of Peter May and
Haynes Miller.  The author is grateful to Bill Dwyer, Igor Kriz,
Gaunce Lewis, and Dennis Sullivan for useful remarks and to
M.~Basterra, P. G. Goerss, D. C. Isaksen, B. Shipley, and J. Wolbert
for helpful comments and suggestions.

%%%%%%%%%%%%%%%%%%%%%%%%%%%%%%%%%%%%%%%%%%%%%%%%%%%%%%%%%%%%%%%%%%%%%%%%
\section{Outline of the Argument}\label{secoutline}

In this section we outline the proofs of Theorems~A and~B.  The
underlying idea is to take advantage of the ``arithmetic square'' to
reduce questions in terms of the integers to questions in terms of the
rational numbers and the prime fields answered in
\cite{einf}.   Let $\Zhat$
denote the finite completion of the integers, i.e., the product of
$\ZPhat$ over all primes $p$, and let $\Qhat=\Zhat\otimes \Q$.  The
integers are the fiber product of the inclusions of $\Zhat$ and $\Q$ in
$\Qhat$, and the square on the left below is called the arithmetic
square.
\[
\diagram
\Z\rto\dto&\Zhat\dto&
&Y\rdotted|>\tip\ddotted|>\tip&\finite{Y}\dto\\
\Q\rto&\Qhat&
&\rat{Y}\rto&\rat{(\finite{Y})}
\enddiagram
\]
The arithmetic square in homotopy theory is the square on the right.
When $Y$ is a finite type nilpotent space, then $Y$ is equivalent to
the homotopy pullback of the rationalization map $\finite{Y}\to
\rat{(\finite{Y})}$ along the rationalization of the finite completion
map $\rat{Y}\to \rat{(\finite{Y})}$.  Looking at the definition of
$\loopf Y$, we obtain an arithmetic square of sorts for $\loopf Y$
from the arithmetic square for $Y$.

\begin{prop}\label{loopfredef}
When $Y$ is finite type nilpotent,
$\loopf Y$ is weakly equivalent to the homotopy pullback of $\Lambda
\finite{Y}\to \rat{(\finite{Y})}$ along $\rat{Y}\to
\rat{(\finite{Y})}$.
\[
\diagram
\loopf Y\rdotted|>\tip\ddotted|>\tip&\Lambda \finite{Y}\dto\\
\rat{Y}\rto&\rat{(\finite{Y})}
\enddiagram
\]
\end{prop}

We would like to describe $\deradj C^{*}Y$ as an analogous homotopy
pullback.  We begin by reviewing the construction of $\deradj$.  We
understand the category $\Ei$ of \einf algebras to be the category of
$\oE$-algebras for some cofibrant \einf operad $\oE$ (over the
integers).  We recall from \cite[\S 2]{einf}, as improved by
\cite{hcofib}, that $\Ei$ is a closed model category with weak
equivalences the quasi-isomorphisms and fibrations the surjections.
The cofibrations are determined by the weak equivalences and the
fibrations, but a concrete description of them may be found in
\cite[2.4--5]{einf}.  We let $\SSet$ denote the category of simplicial
sets.  The normalized cochain functor with coefficients in $\Z$
naturally takes values in the category $\Ei$, and defines a
contravariant functor $C^{*}\colon \SSet\to \Ei$; more generally, for
any commutative ring $R$, the normalized cochains with coefficients in
$R$ defines a contravariant functor $C^{*}(-;R)\colon \SSet\to \Ei$.
These functors have contravariant right adjoints.

\begin{prop}\label{propdeft}
Let $R$ be a commutative ring, let $A$ be an \einf algebra and let
$T(A;R)$ denote the simplicial set of maps
$\Ei(A,C^{*}(\Delta[\ssdot];R))$, where $\Delta[n]$ denotes the
standard $n$-simplex simplicial set.  Then $T(-;R)$ is a contravariant
functor $\Ei\to \SSet$ and is right adjoint to $C^{*}(-;R)$.
\end{prop}

We omit the proof as it is identical to the proof of \cite[4.2]{einf}.
The unit map $Y\to T(C^{*}(Y;R);R)$ is easy to describe:
It takes the $n$-simplex $\sigma$ of $Y$ to the $n$-simplex
\[
C^{*}(f_{\sigma};R)\colon C^{*}(Y;R)\to C^{*}(\Delta[n];R) 
\]
of $T(C^{*}(Y;R);R)$, where $f_{\sigma}\colon \Delta[n]\to Y$ is the
map of simplicial sets corresponding to $\sigma$ (sending the
non-degenerate $n$-simplex of $\Delta[n]$ to $\sigma$).
For the reader familiar with Appendix~A of \cite{einf} and its
notation, we remark that $T(-;\Z)=U(-;\Z)$.  For more general $R$,
$U(-;R)$ was defined as a functor from the category of \einf
$R$-algebras to simplicial sets, so $T(-;R)$ is not the same as
$U(-;R)$; rather, it satisfies $T(A;R)\iso U(A\otimes R;R)$.

The functor $C^{*}(-;R)$ sends weak equivalences of simplicial sets to
quasi-iso\-morphisms of \einf algebras and sends cofibrations (injections)
of simplicial sets to fibrations (surjections) of \einf
algebras.  This implies that ($C^{*}(-,R)$, $T(-;R)$) is a Quillen
adjoint pair, i.e., satisfies the conditions of
\cite[Theorem~4-3]{quil} or \cite[9.7.(i)]{modcat}
(cf. \cite[9.8]{modcat}) that ensure that the right derived functor of
$T(-;R)$ exists and is adjoint to the derived functor of $C^{*}(-;R)$.
We denote this right derived functor as $\Tder(-;R)$.  We remark that since
$T(-;R)$ is contravariant, $\Tder(A;R)$ is constructed by choosing a
\emph{cofibrant} approximation $A'\to A$, and setting
$\Tder(A;R)=T(A';R)$.  As a particular case, we get the functor
$\deradj$ as another name for $\Tder(-;\Z)$.  We summarize this as follows.

\begin{prop}\label{propderadj}
The functor $T(-;R)$ converts cofibrations of \einf algebras to Kan
fibrations of simplicial sets and converts quasi-isomorphisms of
cofibrant \einf algebras to weak equivalences of simplicial sets.
The right derived functor $\Tder(-;R)$ of $T(-;R)$ exists and is right
adjoint to $C^{*}(-;R)$.
\end{prop}

\begin{notn}
We write $\adj$ for $T(-;\Z)$ and $\deradj$ for $\Tder(-;\Z)$.
\end{notn}

The following theorem proved in Section~\ref{sechocar} gives the
arithmetic square we need.

\begin{thm}\label{arithsquare}
Let $A$ be a cofibrant \einf algebra.  Then the diagram
\[ 
\diagram
U(A)\rto\dto&T(A;\Zhat)\dto\\
T(A;\Q)\rto&T(A;\Qhat)
\enddiagram
\]
is homotopy cartesian, i.e., the induced map from $U(A)$ to the
homotopy pullback of $T(A;\Zhat)\to T(A;\Qhat)$ along $T(A;\Q)\to
T(A;\Qhat)$ is a weak equivalence.
\end{thm}

The proof of Theorem~B involves showing that when we take $A$ to be a
cofibrant approximation of $C^{*}Y$, the square in
Theorem~\ref{arithsquare} is equivalent to the square in Proposition~\ref{loopfredef} describing
$\loopf Y$.  The rational part is
straightforward to analyze.  The adjunction of
Proposition~\ref{propdeft} for $U=T(-;\Z)$ gives us a map
$Y\to U(C^{*}Y)$, which
we compose to obtain a map
\begin{equation}\label{uniteta}
\eta \colon Y\to U(C^{*}Y)=T(C^{*}(Y);\Z)\to T(C^{*}Y;R) \to T(A;R).
\end{equation}
Passing to the homotopy category, we obtain a map $Y\to
\Tder(C^{*}Y;R)$ that is natural in $Y$.  
The following theorem is essentially \cite[A.6]{einf}; we review the
details and give a complete argument in Section~\ref{secpostnikov}. 

\begin{thm}\label{rationalpieceone}
Let $Y$ be a connected finite type nilpotent simplicial set.  Then the
natural map in the homotopy category $\eta \colon Y\to \Tder(C^{*}Y;\Q)$ is
rationalization. 
\end{thm}

To analyze $\Tder(C^{*}Y;\Zhat)$, we use the following theorem proved
in Sections~\ref{secwittring} and~\ref{secpostnikov}.

\begin{thm}\label{wittring}
Let $R$ be a complete discrete valuation ring with maximal ideal
$\mathfrak{m}$, such that the residue field $R/\mathfrak{m}$ has
finite characteristic. Then the natural map
\[
\Tder(C^{*}Y;R) \to \Tder(C^{*}Y;R/\mathfrak{m})
\]
is a weak equivalence when $Y$ is connected finite type nilpotent.
\end{thm}

As a consequence, the map $\Tder(C^{*}Y;\ZPhat)\to \Tder(C^{*}Y;\FP)$
is a weak equivalence.  Appendix~A of \cite{einf} proves (in our
current notation) that when $Y$ is connected finite type nilpotent,
$\Tder(C^{*}Y;\FP)$ is weakly equivalent to $\Lambda \finite{Y}_{p}$,
the free loop space of the $p$-completion of $Y$.  Since $\Zhat$ is
the product of $\ZPhat$ over all primes $p$, and $\Lambda
\finite{Y}$ is equivalent to the product of $\Lambda \finite{Y}_{p}$
over all primes $p$, we obtain the following corollary.

\begin{cor}\label{loopcor}
Let $Y$ be a connected finite type nilpotent simplicial set.  Then
there is natural isomorphism in the homotopy category $\Lambda
\finite{Y}\to \Tder(C^{*}Y;\Zhat)$.
\end{cor}

In order to describe the last corner, we need an observation on the 
natural map $\eta$ when $R$ contains $\Zhat$.
In the description of $\eta$ in (\ref{uniteta}), we used the unit of
the ($C^{*}$, $U$) adjunction; however, it is clear from the explicit
description of the unit map above that when the commutative ring $R$
contains $\Zhat$, we obtain the same map as the following composite
using the ($C^{*}(-;\Zhat)$, $T(-;\Zhat)$) adjunction.
\[
Y\to T(C^{*}(Y;\Zhat);\Zhat)\to T(C^{*}Y;R) \to T(A;R)
\]
Factoring the map $A\to C^{*}(Y;\Zhat)$ as a cofibration $A\to
\finite{A}$ followed by an acyclic fibration $\finite{A}\to
C^{*}(Y;\Zhat)$, for some \einf algebra $\finite{A}$, we see that the
map $\eta$ factors through the map $Y\to
T(C^{*}(\finite{A};\Zhat);\Zhat)$. In other 
words, $\eta$ may be written as a composite
\[
Y\to \Tder(C^{*}(Y;\Zhat);\Zhat) \to \Tder(C^{*}Y;R),
\]
where the first map is the unit of the derived adjunction
($C^{*}(-;\Zhat)$, $\Tder(-;\Zhat)$).
The significance of this is that the finite completion map $Y\to
\finite{Y}$ induces a quasi-isomorphism $C^{*}(\finite{Y};\Zhat)\to
C^{*}(Y;\Zhat)$ when $Y$ is finite type nilpotent (or when we take
$\finite{Y}$ to denote Bousfield finite completion \cite{bousloc}).
It follows that the induced map $\Tder(C^{*}(Y;\Zhat);\Zhat)\to
\Tder(C^{*}(\finite{Y};\Zhat);\Zhat)$ is a weak equivalence, and hence
that the map $\eta$ factors in the homotopy category as the
composite of finite completion $Y\to \finite{Y}$ and the natural map
\begin{equation}\label{etahat}
\finite{\eta}\colon \finite{Y}\to
\Tder(C^{*}(\finite{Y};\Zhat);R)\homeq
\Tder(C^{*}(Y;\Zhat);R)\to
\Tder(C^{*}Y;R)
\end{equation}
The following theorem regarding this map is proved in
Section~\ref{secpostnikov}.

\begin{thm}\label{rationalpiecetwo}
Let $Y$ be a connected finite type nilpotent simplicial set.
Then the natural map in the homotopy category $\finite{\eta}\colon
\finite{Y}\to \Tder(C^{*}Y;\Qhat)$ is rationalization.
\end{thm}

We note that the universal property of rationalization then implies
that the following diagram commutes.
\[
\diagram
Y\rto\dto&\rat{Y}\rto^(.3){\sim}\dto&\Tder(C^{*}Y;\Q)\dto\\
\relax \finite{Y}\rto&\rat{\finite{Y}}\rto_(.3){\sim}&\Tder(C^{*}Y;\Qhat)
\enddiagram
\]
This identifies one map in the diagram in Theorem~\ref{arithsquare}.
Identifying the other map in the diagram is significantly more difficult.  The
following theorem is proved in Section~\ref{secyoneda}.

\begin{thm}\label{looppiece}
Let $Y$ be a connected finite type nilpotent simplicial set.  The
natural isomorphism in the homotopy category $\Lambda \finite{Y}\to
\Tder(C^{*}Y;\Zhat)$ makes the following diagram in the homotopy
category commute,
\[
\diagram
\relax \Lambda \finite{Y}\dto\rrto^{\sim}&&\Tder(C^{*}Y;\Zhat)\dto\\
\relax \finite{Y}\rto&\rat{\finite{Y}}\rto_(.3){\sim}&\Tder(C^{*}Y;\Qhat)
\enddiagram
\]
\end{thm}

Theorem~B, which requires no naturality, is an easy consequence of
Theorems~\ref{arithsquare}, \ref{rationalpieceone},
\ref{rationalpiecetwo}, and~\ref{looppiece}
in the connected case; the non-connected case follows from the
connected case and \cite[3.1]{finite} (for $k=\Z$).  If we compose the
inverse in the homotopy category of the equivalence $\loopf Y\to
\deradj C^{*}Y$ with the map $\loopf Y\to Y$, we obtain a map in the
homotopy category $\deradj (C^{*}Y)\to Y$.  The preceding theorems
are not strong enough to imply that we can arrange for this map to be
natural in $Y$; however, the proof of these theorems is.  We
prove the following theorem in Section~\ref{secthmai}.

\begin{thm}\label{thmai}
There is a natural map in the homotopy category $\epsilon \colon
\deradj(C^{*}Y)\to Y$ for $Y$ a finite type nilpotent simplicial set,
such that the composite with the unit 
\[
Y \overto{\hphantom{\eta}} \deradj(C^{*}Y) \overto{\epsilon} Y
\]
is the identity.
\end{thm}

This theorem gives the first statement of Theorem~A, the natural
retraction $\retr\colon \HoEiMap(C^{*}Y,C^{*}X)\to \HoSPMap(X,Y)$.  The
remainder of Theorem~A requires us to show that maps $f,g\in
\HoEiMap(C^{*}Y,C^{*}X)$ that satisfy $\retr(f)=\retr(g)$ must satisfy
$H^{*}f=H^{*}g$.  We do this by identifying in \einf algebra terms
when maps satisfy $\retr(f)=\retr(g)$.  We need a ring that plays the
same role for $\ZPhat$ that $\FPbar$ plays for $\FP$.  This ring is the
Witt vectors of $\FPbar$.  For a proof of the uniqueness statement in
the following definition, see \cite[II\S5]{serrelocalfields}; for an
elementary construction of the Witt vectors, see
\cite[II\S6]{serrelocalfields}.

\begin{defn}\label{defw}
Let $W(\FPbar)$ denote the $p$-typical Witt vectors of $\FPbar$, the
unique complete discrete valuation ring with maximal ideal $(p)\neq
(0)$ and residue field~$\FPbar$.
%Let $\Phi_{p}$ denote the Frobenius
%automorphism of $W(\FPbar)$, the unique endomorphism that induces the
%Frobenius on the residue field.  
Let $\W = \prod_{p} W(\FPbar)$. %, and
%let $\Phi = \prod_{p} \Phi_{p}\colon \W\to \W$.
\end{defn}

By Theorem~\ref{wittring} and the Main Theorem of \cite{einf}, we have
that the natural map $Y\to \Tder(C^{*}Y;W(\FPbar))$ is $p$-completion
when $Y$ is connected finite type nilpotent.  It follows that the
natural map $Y\to \Tder(C^{*}Y;\W)$ is finite completion when $Y$ is
connected finite type nilpotent.  Just as in \cite{einf}, we may use
the closure of $\FP$ under degree $p$ extensions in place of $\FPbar$,
i.e., the (infinite) algebraic extension that is the fixed field of
$\myop\prod_{\ell\neq p}\Zhat_{\ell} < \Zhat = \Gal(\FPbar/\FP)$.

In Section~\ref{secyoneda}, we prove the following theorem.

\begin{thm}\label{thma0}
Let $Y$ be a connected finite type nilpotent simplicial set.  There is
an isomorphism in the homotopy category $\finite{Y}\to T(C^{*}Y;\W)$
such that the following diagram in the homotopy category commutes
\[
\diagram
\relax \Lambda \finite{Y}\rto^{\sim}\dto&\Tder(C^{*}Y;\Zhat)\dto\\
\relax \finite{Y}\rto_{\sim}&\Tder(C^{*}Y;\W)
\enddiagram
\]
where the top map is the map in Theorem~\ref{looppiece}.
\end{thm}

Given a map $f\colon C^{*}Y\to C^{*}X$, write $\Wmap{f}$ and
$\Qmap{f}$ for the composite maps 
\[
C^{*}Y\overto{f}C^{*}X\to C^{*}(X;\W),\qquad 
C^{*}Y\overto{f}C^{*}X\to C^{*}(X;\Q),
\]
respectively.  We have the following corollary of the previous theorem.

\begin{cor}\label{coraii}
Let $X$ be a simplicial set, let $Y$ be a finite type
nilpotent simplicial set, and let $f$ and $g$ be maps in $\HoEiMap(C^{*}Y, 
C^{*}X)$ such that $\retr(f)=\retr(g)$ in
$\HoSPMap(X,Y)$. Then $\Wmap{f}=\Wmap{g}$ in
$\HoEiMap(C^{*}Y,{C^{*}(X;\W)})$ and $\Qmap{f}=\Qmap{g}$ in
$\HoEiMap(C^{*}Y,{C^{*}(X;\Q)})$.
\end{cor}

\begin{proof}
If $\retr(f)=\retr(g)$,  then $f$ and $g$ are sent
to the same element under the map
\[
\HoEiMap(C^{*}Y,C^{*}X)\iso \HoSPMap(X,\deradj C^{*}Y)\to
\HoSPMap(X,\finite{Y}) \times \HoSPMap(X,\rat{Y}).
\]
Since we can decompose the above map into a product over the
components of $X$, it suffices to consider the case when $X$ is
connected.  We can then decompose into a disjoint union
over the components of $Y$ (by \cite[3.1]{finite}), and so it suffices
to consider the case when $Y$ is also connected.

By the previous theorem, $f$ and $g$ are sent to the same element
under the map above if and only if $f$ and $g$ are
sent to the same element under the map
\[
\HoEiMap(C^{*}Y,C^{*}X)\iso \HoSPMap(X,\deradj C^{*}Y)\to
\HoSPMap(X,{\Tder(C^{*}Y;\W)}) \times \HoSPMap(X,{\Tder(C^{*}Y;\Q)}).
\]
By the adjunction of Proposition~\ref{propderadj}, this is equivalent
to the condition that $\Wmap{f}=\Wmap{g}$ and $\Qmap{f}=\Qmap{g}$.
\end{proof}

We get the second statement of Theorem~B as a consequence: When $X$
is finite type, the map $C^{*}X\to C^{*}(X;\W)\times C^{*}(X,\Q)$ is
injective on cohomology.

%%%%%%%%%%%%%%%%%%%%%%%%%%%%%%%%%%%%%%%%%%%%%%%%%%%%%%%%%%%%%%%%%%%%%%%%
\section{Simplicial and Cosimplicial Resolutions}\label{sechocar}

In this section we prove Theorem~\ref{arithsquare}, which establishes
the arithmetic square fracturing the functor $\deradj$.  The proof is
a standard argument using the tools of simplicial and cosimplicial
resolutions introduced in \cite{dk3}.  We review the basic definitions 
and terminology, which we use throughout the remainder of the paper.

Let $\Modcat$ be a closed model category.
For an object $X$ of $\Modcat$, a \term{cosimplicial resolution of
$X$} is a cosimplicial object $X\supdot$ together with a
weak equivalence $X^{0}\to X$ such that $X^{0}$ is cofibrant, each
coface map in $X\supdot$ is an acyclic cofibration, and each map
$L^{n}\to X^{n+1}$ is a cofibration, where $L^{n}$ is
the object denoted as $(d^{*},X^{n})$ in 
\cite[4.3]{dk3}: It is defined to be the colimit of the diagram in 
$\Modcat$ with objects
\begin{itemize}
\item For each $i$, $0\leq i\leq n+1$, a copy of $X^{n}$ labeled
$(d^{i},X^{n})$ 
\item For each $(i,j)$, $0\leq i<j\leq n+1$, a copy of $X^{n-1}$ labeled
$(d^{j}d^{i},X^{n-1})$ (we understand $X^{-1}$ to be the initial object). 
\end{itemize}
and maps
\begin{itemize}
\item For each $(i,j)$, $0\leq i<j\leq n+1$, a map
$(d^{j}d^{i},X^{n-1})\to (d^{j},X^{n})$ given by the map $d^{i}\colon
X^{n-1}\to X^{n}$.
\item For each $(i,j)$, $0\leq i<j\leq n+1$, a map
$(d^{j}d^{i},X^{n-1})\to (d^{i},X^{n})$ given by the map $d^{j-1}\colon
X^{n-1}\to X^{n}$.
\end{itemize}
Simplicial resolutions are defined dually; a simplicial resolution
in $\Modcat$ is a cosimplicial resolution in $\Modcat^{\op}$.
The following proposition and its proof give enlightening examples of
simplicial and cosimplicial resolutions.

\begin{prop}\label{cstarsimpres}
For any commutative ring $R$, $C^{*}(\Delta[\ssdot];R)$ is a
simplicial resolution of $R$ in the category of \einf algebras.
\end{prop}

\begin{proof}
Clearly, $\Delta[\ssdot]$ is a cosimplicial resolution of $*$ in the
category of simplicial sets: The object $L^{n}$ is just the boundary
of $\Delta[n+1]$.  If we view $C^{*}(-;R)$ as a functor to
$\Ei^{\op}$, then it is a Quillen left adjoint.  In particular it
preserves cofibrations, acyclic cofibrations, and colimits of finite
diagrams.  It follows that $C^{*}(-;R)$ takes cosimplicial resolutions
of simplicial sets to cosimplicial resolutions in $\Ei^{\op}$.  
\end{proof}

Cosimplicial resolutions have the following basic properties.

\begin{prop}\label{preserveswe}
Let $X\supdot$ be a cosimplicial resolution. The functor
$\Modcat(X\supdot,-)$ from $\Modcat$ to simplicial sets
preserves fibrations and acyclic fibrations, and preserves weak
equivalences between fibrant objects.
\end{prop}

\begin{proof}
If we write $L^{n,i}$ for the colimit of the diagram analogous to the
one above but omitting $(d^{i},X^{n})$, a straightforward induction
(written in detail in \cite[\S 6]{hh}) shows that the map $L^{n,i}\to
X^{n+1}$ is an acyclic cofibration.  When $f\colon Y\to Z$ is a
fibration, the left lifting property for the acyclic cofibrations
$L^{n,i}\to X^{n+1}$ with respect to the fibration $f$ translates
under the universal property of the colimits defining the $L^{n,i}$ into
the extension condition for $\Modcat(X\supdot,f)$ to be a Kan
fibration.  Likewise, when $f$ is an acyclic fibration, the left
lifting property for the cofibration $L^{n}\to X^{n+1}$ with respect
to $f$ translates under the universal property of the
colimit defining $L^{n}$ into the extension condition for
$\Modcat(X\supdot,f)$ to be an acyclic Kan fibration.
K.~Brown's lemma \cite[9.9]{modcat} then implies that
$\Modcat(X\supdot,-)$ preserves weak equivalences between fibrant
objects. 
\end{proof}

\begin{prop}\label{compareresolutions}
Let $X\supdot$ be a cosimplicial resolution of $X$ and let $Y\subdot$ be a
simplicial resolution of $Y$.  If $X$ is cofibrant and $Y$ is fibrant,
then the maps of simplicial sets
\[
\Modcat(X\supdot,Y)\to 
\diag \Modcat(X\supdot,Y\subdot)\from
\Modcat(X,Y\subdot)
\]
are weak equivalences.
\end{prop}

\begin{proof}
When we regard $\Modcat(X\supdot,Y)$ as a bisimplicial set, constant
in the second simplicial direction, then the map of bisimplicial sets
$\Modcat(X\supdot,Y)\to \Modcat(X\supdot,Y\subdot)$ has the property
that each map $\Modcat(X\supdot,Y)\to \Modcat(X\supdot,Y_{n})$ is a
weak equivalence by the previous proposition.  It follows that the map
$\Modcat(X\supdot,Y)\to\diag\Modcat(X\supdot,Y\subdot)$ is a weak
equivalence.  The map
$\Modcat(X,Y\subdot)\to\Modcat(X\supdot,Y\subdot)$ is a weak
equivalence by the same argument in $\Modcat^{\op}$.
\end{proof}

\begin{prop}
Every object has a simplicial resolution.
\end{prop}

\begin{proof}
This is a straightforward factorization argument; see \cite[6.7]{dk3}
for details.
\end{proof}

For the application of cosimplicial resolutions to the proof of
Theorem~\ref{arithsquare}, it is convenient to introduce the following
terminology.  Let 
\begin{equation}\label{hocasq}
{\diagram
W\rto\dto&X\dto\\Y\rto&Z
\enddiagram}
\end{equation}
be a commutative square in a closed model category $\Modcat$. 
By factoring maps by weak equivalences followed by fibrations, we can
form a commutative square 
\[
{\diagram
X\rto\dto_{\sim}&Z\dto_{\sim}&Y\dto_{\sim}\lto\\
X'\rfib&Z'&Y'\lfib
\enddiagram}
\]
where $Z'$ is fibrant, $X'\to Z'$ and $Y'\to Z'$ are fibrations, and
the maps $X\to X'$, $Y\to Y'$, and $Z\to Z'$ are weak equivalences.
We have an induced map from $W$ to the pullback $X'\times_{Z'}Y'$.  
If this map is a weak equivalence for some choice of $X'\to Z'\from
Y'$, it is a weak equivalence for any choice.

\begin{defn}\label{defhomotopycartesian}
We say that the square (\ref{hocasq}) is \term{homotopy cartesian} if
the map $W\to X'\times_{Z'}Y'$ is a weak equivalence.
\end{defn}

This is equivalent to the usual definition in the category of spaces
or simplicial sets.  Maps out of cosimplicial resolutions preserve the
homotopy cartesian property when the objects in the square are fibrant. 

\begin{prop}\label{homotopycartesian}
Let $A\supdot$ be a cosimplicial resolution.  If the square on the
left is homotopy cartesian in $\Modcat$, 
\[
\diagram
W\rto\dto&X\dto&
&\Modcat(A\supdot,W)\dto\rto&\Modcat(A\supdot,X)\dto\\
Y\rto&Z&
&\Modcat(A\supdot,Y)\rto&\Modcat(A\supdot,Z)
\enddiagram
\]
and $W$, $X$, $Y$, $Z$ are all
fibrant, then the square on the right is homotopy cartesian in
$\SSet$. 
\end{prop}

\begin{proof}
Let $X'\to Z'\from Y'$ be as in Definition~\ref{defhomotopycartesian}
above, and let $W'$ be the pullback $X'\times_{Z'}Y'$.  Then the
square 
\[
\diagram
\Modcat(A\supdot,W')\dfib\rfib&\Modcat(A\supdot,X')\dfib\\
\Modcat(A\supdot,Y')\rfib&\Modcat(A\supdot,Z')
\enddiagram
\]
is a pullback square with all maps fibrations and all objects Kan
complexes.  The vertical maps
\[
\diagram
\Modcat(A\supdot,X)\rto\dto_{\sim}
&\Modcat(A\supdot,Z)\dto_{\sim}
&\Modcat(A\supdot,Y)\lto\dto_{\sim}\\
\Modcat(A\supdot,X')\rfib
&\Modcat(A\supdot,Z')
&\Modcat(A\supdot,Y')\lfib
\enddiagram
\]
are weak equivalences, as is the map $\Modcat(A\supdot,W)\to
\Modcat(A\supdot,W')$.  This proves that the square on the right in the statement
is homotopy cartesian.
\end{proof}

The following proposition is the last fact we need for
the proof of Theorem~\ref{arithsquare}.

\begin{prop}\label{ashc}
The arithmetic square of $\Z$, $\Zhat$, $\Q$, $\Qhat$ is
homotopy cartesian in the category of \einf algebras.
\end{prop}

\begin{proof}
Factor $\Zhat\to \Qhat$ and $\Qhat$ through acyclic cofibrations
$\Zhat\to B$ and $\Qhat\to C$ followed by fibrations $B\to \Qhat$ and
$C\to \Qhat$, and let $D$ be the pullback $B\times_{\Qhat}C$.  
We get an induced long exact sequence on cohomology 
\[
\cdots \to H^{-1}(\Qhat)\to H^{0}(D) \to (H^{0}B \times H^{0}C)
\to H^{0}(\Qhat) \to H^{1}(D) \to \cdots.
\]
Of course $H^{-1}\Qhat=0$ and the map 
\[
\Zhat \times \Q = H^{0}B \times H^{0}C
\to  H^{0}(\Qhat) = \Qhat
\]
is surjective, so $H^{1}D=0$.  The long exact sequence above is
therefore isomorphic to the lengthened short exact sequence 
\[
\cdots \to 0 \to \Z \to (\Zhat\times \Q)\to \Qhat\to 0\to \cdots 
\]
and the map $\Z\to D$ is a quasi-isomorphism.
\end{proof}

We can now prove Theorem~\ref{arithsquare}.

\begin{proof}[Proof of Theorem~\ref{arithsquare}]
We choose a cosimplicial resolution $A\supdot$ of $A$.
Then Propositions~\ref{ashc} and~\ref{homotopycartesian} imply that
the square 
\[
\diagram
\Ei(A\supdot,\Z)\rto\dto&\Ei(A\supdot,\Zhat)\dto\\
\Ei(A\supdot,\Q)\rto&\Ei(A\supdot,\Qhat)
\enddiagram
\]
is homotopy cartesian.  The theorem now follows from
Propositions~\ref{cstarsimpres} and~\ref{compareresolutions}. 
\end{proof}

%%%%%%%%%%%%%%%%%%%%%%%%%%%%%%%%%%%%%%%%%%%%%%%%%%%%%%%%%%%%%%%%%%%%%%%%
\section{A Reduction of Theorem~\ref{wittring}}
\label{secwittring}

In this section, we use the theory of resolutions reviewed in the
previous section to reduce Theorem~\ref{wittring} to a statement only
involving the residue field $k=R/\mathfrak{m}$.  
The basic idea is that a complete
discrete valuation ring is the limit of a sequence of ``square zero
extensions'' by the quotient field and that such extensions are
obtained by base change from a trivial extension, denoted $k\oplus
k[-1]$ below.  As we explain below, this reduces the problem to
understanding the simplicial sets of maps into $k\oplus k[-1]$
factoring a given map into $k$.  The main reduction,
Theorem~\ref{aqvanish} below, states that these are contractible; we
prove this theorem in the next section.

Because of the structure of the argument that follows, it is convenient to
reformulate Theorem~\ref{wittring} in the following form, where the
basic construction is functorial in \einf algebra maps of $R$ instead
of merely commutative ring maps of $R$; it is equivalent to the
original statement by Proposition~\ref{compareresolutions}.

\begin{thm}\label{wittringtwo}
Let $Y$ be a connected finite type nilpotent simplicial set and let
$A\supdot$ be a cosimplicial resolution of $C^{*}Y$.  Let
$(R,\mathfrak{m})$ be a complete discrete valuation ring whose residue
field $k=R/\mathfrak{m}$ has finite characteristic.  Then the map of
simplicial sets $\Ei(A\supdot,R)\to \Ei(A\supdot,k)$ is a
weak equivalence. 
\end{thm}

With notation as above, choose $\pi$ to be an irreducible element of
$\mathfrak{m}$; then $\mathfrak{m}=(\pi)$. Since the valuation is
complete, the canonical map 
\[
R \to \Lim R/(\pi^{n})
\]
is an isomorphism. Since the maps $R/(\pi^{n+1})\to R/(\pi^{n})$ are
surjective, we have 
\[
\Ei(A\supdot,R)\iso \Lim \Ei(A\supdot,R/(\pi^{n}))
\]
is the limit of Kan fibrations of Kan complexes.  Thus, to prove
Theorem~\ref{wittringtwo}, it suffices to show that each map
\[
\Ei(A\supdot,R/(\pi^{n+1}))\to \Ei(A\supdot,R/(\pi^{n}))
\]
is a weak equivalence.

Consider the following variant of the tower $R/(\pi^{n})$.  Let
$R_{n}$ be the Koszul complex associated to $\pi^{n}$:  This is the
commutative differential graded $R$-algebra that is the exterior
$R$-algebra on an element $x_{n}$ whose differential is $\pi^{n}$.
We have a map of commutative differential graded $R$-algebras
$R_{n+1}\to R_{n}$ obtained by sending $x_{n+1}$ to $\pi\cdot x_{n}$;
this induces on homology the map of commutative rings
$R/(\pi^{n+1})\to R/(\pi^{n})$ 

Let $k\oplus k[-1]$ denote the graded commutative $k$-algebra, exterior
on an element of degree $-1$ (where differentials raise degree).  We
regard this as a commutative differential graded $R$-algebra with zero
differential, augmented to $k$ by the projection map sending the
exterior generator to zero.  (We think of $k\oplus k[-1]$ as the trivial
``square zero extension'' of $k$.)  We have a map of commutative
differential graded $R$-algebras $R_{n}\to k\oplus k[-1]$ that sends
$x_{n}$ to the exterior generator; this map is surjective.

\begin{prop}
The commutative diagram 
\[
\diagram
R_{n+1}\rfib\dto&k\dto\\R_{n}\rfib&k\oplus k[-1]
\enddiagram
\]
is a pullback square and homotopy cartesian square in the category of
\einf algebras.
\end{prop}

By Proposition~\ref{homotopycartesian}, it follows that the square
\[
\diagram
\Ei(A\supdot,R_{n+1})\rto\dto&\Ei(A\supdot,k)\dto\\
\Ei(A\supdot,R_{n})\rto&\Ei(A\supdot,k\oplus k[-1])\\
\enddiagram
\]
is homotopy cartesian in the category of simplicial sets.
Since the maps
\[
\Ei(A\supdot,R_{n+1})\to \Ei(A\supdot,R/(\pi^{n+1})), \qquad 
\Ei(A\supdot,R_{n})\to \Ei(A\supdot,R/(\pi^{n}))
\]
are weak equivalences, we are reduced to showing that the map 
\[
\Ei(A\supdot,k)\to \Ei(A\supdot,k\oplus k[-1])
\]
is a weak
equivalence, or equivalently, that the retraction
\[
r\colon \Ei(A\supdot,k\oplus k[-1])\to \Ei(A\supdot,k)
\]
(induced by the augmentation $k\oplus k[-1]\to k$)
is a weak equivalence.

Since the map $k\oplus k[-1]\to k$ is surjective, the map $r$ is a
Kan fibration.  It therefore suffices to show that the fiber of $r$ at
each 
vertex is contractible.  Choose and fix a map $b\colon A^{0}\to k$.  If we
write $\Ei/k$ for the category of \einf algebras lying over $k$, then 
the fiber of $r$ over the point $b$ is exactly the simplicial set
\[
(\Ei/k)(A\supdot,k\oplus k[-1]),
\]
where the map $A^{s}\to k$ is the composite of the degeneracy $A^{s}\to
A^{0}$ and the map $b\colon A^{0}\to k$.  Thus, Theorem~\ref{wittring}
reduces to showing that this simplicial set is contractible.  We
reformulate this as Theorem~\ref{aqvanish} below, using
Proposition~\ref{compareresolutions} to switch from using a
cosimplicial resolution of $C^{*}Y$ to using a simplicial resolution
of $k\oplus k[-1]$.  

We construct a simplicial resolution for $k\oplus k[-1]$ as follows.
For a simplicial set $X$, let $C^{*}(X;k[-1])$ denote the differential
graded $k$-module obtained from $C^{*}(X;k)$ by shifting one degree down,
\[
C^{n}(X;k[-1]) = C^{n+1}(X;k).
\]
We make $k\oplus C^{*}(X;k[-1])$ an \einf algebra lying over $k$ by
giving it the \term{square zero multiplication}.  This is an \einf
algebra structure coming from the augmented commutative differential
graded $k$-algebra structure where any pair of elements in the
augmentation ideal $C^{*}(X;k[-1])$ multiply to zero
(i.e., $C^{*}(X;k[-1])$ is a square zero ideal).  The elements in $k$
multiply normally and multiply with elements of $C^{*}(X;k[-1])$ by
the usual $k$-module action.  Now consider the simplicial object of
$\Ei/k$ given by
\[
k\oplus C^{*}(\Delta[\ssdot];k[-1]);
\]
clearly, this is a simplicial resolution of $k\oplus k[-1]$ in $\Ei/k$.

\begin{defn}
Let $A$ be an \einf algebra lying over $k$.  Let 
\[
\AQS(A;k) = (\Ei/k)(A,k\oplus C^{*}(\Delta[\ssdot];k[-1]))
\]
\end{defn}

The work above together with Proposition~\ref{compareresolutions}
reduces Theorem~\ref{wittring} to the following theorem.

\begin{thm}\label{aqvanish}
Let $k$ be a field of positive characteristic.
Let $Y$ be a connected finite type nilpotent simplicial set, let
$C\to C^{*}Y$ be a cofibrant approximation, and let $C\to k$ be a map
of \einf algebras.  Then $\AQS(C;k)$ is contractible.
\end{thm}

As stated, the theorem above requires us to prove contractibility for
all cofibrant approximations $C$.  The following lemma allows us to
replace the implicit ``for every'' with a ``there exists'' and work
with a convenient cofibrant approximation instead of an arbitrary one.
It is an immediate consequence of \cite[2.12]{einf} or the fact that
in the model category $\Ei$ all objects are fibrant.

\begin{lem}\label{maptoc}
Let $C\to C^{*}Y$ be a given cofibrant approximation.  If $A\to
C^{*}Y$ is any cofibrant approximation, then there exists a
quasi-isomorphism $A\to C$. 
\end{lem}

%%%%%%%%%%%%%%%%%%%%%%%%%%%%%%%%%%%%%%%%%%%%%%%%%%%%%%%%%%%%%%%%%%%%%%%%
\section{The Proof of Theorems~\ref{rationalpieceone}, 
\ref{rationalpiecetwo}, and~\ref{aqvanish}}\label{secpostnikov}

In this section we prove Theorems~\ref{rationalpieceone}
and~\ref{rationalpiecetwo} from Section~\ref{secoutline}, and
Theorem~\ref{aqvanish} from the previous section.  What these theorems
have in common is that they are proved by induction up a principally
refined Postnikov tower of a connected finite type nilpotent
simplicial set. 

Recall that a principally refined Postnikov tower is a finite or
infinite tower of fibrations
\[
\cdots \to Y_{i}\to \cdots \to Y_{1}\to Y_{0}=*
\]
such that each $Y_{i}$ is formed as the pullback of a map $Y_{i}\to
K(G_{i},n_{i}+1)$ along the fibration $L(G_{i},n_{i}+1)\to
K(G_{i},n_{i}+1)$, and $n_{1}$, $n_{2}$, \dots is a non-decreasing
sequence of positive integers, taking on a given value at most
finitely many times.  Here $K(G_{i},n_{i}+1)$ denotes the (standard)
Eilenberg--Mac\,Lane complex with $\pi_{n_{i}+1}=G_{i}$, the map
$L(G_{i},n_{i}+1)\to K(G_{i},n_{i}+1)$ is a Kan fibration, and
$L(G_{i},n_{i}+1)$ is a contractible Kan complex.  See for example,
\cite[\S 23]{simpob}.

A connected simplicial set $Y$ is finite type nilpotent if and only if
there exists a principally refined Postnikov tower $\{Y_{i}\to
K(G_{i},n_{i}+1) \}$ and a weak equivalence $Y\to\Lim Y_{i}$ such that
each $G_{i}$ is $\Z$ or $\ZP$ for some $p$ (depending on $i$).  If
we let $\finite{G}_{i}$ denote the pro-finite completion of $G_{i}$,
then $\finite{Y}$, the finite completion of $Y$, admits a principally
refined Postnikov tower $\{\finite{Y}_{i}\to
K(\finite{G}_{i},n_{i}+1)\}$. Specifically, each $\finite{Y}_{i}$ is
the finite completion of $Y_{i}$ and the maps $\finite{Y}_{i}\to
K(\finite{G}_{i},n_{i}+1)$ represent the finite completion of the maps
$Y_{i}\to K(G_{i},n_{i}+1)$, i.e., the diagram
\[
\diagram
Y_{i}\rto\dto&K(G_{i},n_{i}+1)\dto\\
\finite{Y}_{i}\rto&K(\finite{G}_{i},n_{i}+1)
\enddiagram
\]
commutes in the homotopy category.  The following proposition gives
the first reduction. 

\begin{prop}\label{towerprop}
Let $\cdots \to Y_{1}\to Y_{0}$ be a principally refined Postnikov
tower for the connected finite type 
nilpotent simplicial set $Y$.
If each $Y_{i}$ satisfies Theorem~\ref{rationalpieceone},
\ref{rationalpiecetwo}, or~\ref{aqvanish}, then so does $Y$.
\end{prop}

\begin{proof}
Let $A_{0}\to C^{*}Y_{0}$ be a cofibrant approximation, and
inductively construct $A_{i}\to C^{*}Y_{i}$ by factoring the composite
map $A_{i-1}\to C^{*}Y_{i-1}\to C^{*}Y_{i}$ as a cofibration
$A_{i-1}\to A_{i}$ followed by an acyclic fibration $A_{i}\to
C^{*}Y_{i}$.  Let $A=\Colim A_{i}$; then the canonical map $A\to
C^{*}(\Lim Y)\to C^{*}Y$ is a quasi-isomorphism since the canonical
map $\Colim H^{*}(\Lim Y_{i})\to H^{*}Y$ is an isomorphism.  Let
\[
T_{i}=T(A_{i},\Q),\qquad
\finite{T}_{i}=T(A_{i},\Qhat),\qquad
\AQS_{i}=\AQS(A_{i},k),
\]
and let $T$, $\finite{T}$, and $\AQS$ be the corresponding
constructions for $A$.  In the last case we choose the augmentations
for the $A_{i}\to k$ by choosing an augmentation for $A\to k$ using
Lemma~\ref{maptoc} to find a quasi-isomorphism $A\to C$ from $A$ to
our given cofibrant approximation $C$ and composing with the given
augmentation $C\to k$. 
Then we have 
\[
T\iso \Lim T_{i},\qquad 
\finite{T}\iso \Lim \finite{T}_{i},\qquad 
\AQS\iso \Lim \AQS_{i},
\]
the limits of towers of fibrations of Kan complexes.  We have a
commuting diagram of towers $Y_{i}\to T_{i}$, and the natural maps in
the homotopy category $\finite{\eta}$ give us a homotopy commuting
diagram of towers $\finite{Y}_{i}\to \finite{T}_{i}$, where
$\{\finite{Y}_{i}\}$ is the Postnikov tower for
$\finite{Y}$ corresponding to $\{ Y_{i}\}$.  By hypothesis the $T_{i}$,
$\finite{T}_{i}$, or $\AQS_{i}$ are all connected and the sequence
$\pi_{n}T_{i}$, $\pi_{n}\finite{T}_{i}$, or $\pi_{n}\AQS_{i}$ is
Mittag-Leffler, and so we have
\[
\pi_{n}T\iso \Lim \pi_{n}T_{i},\qquad 
\pi_{n}\finite{T}\iso \Lim\pi_{n}\finite{T}_{i},\qquad 
\pi_{n}\AQS\iso \Lim \pi_{n}\AQS_{i}.
\]
The proposition now follows.
\end{proof}

Next we need the following proposition that explains the effect on
cochain \einf algebras of pullback along a Kan fibration.
It is not formal; it is proved by the same argument as
\cite[5.2]{einf}, which we omit.  In the statement, \term{homotopy
cocartesian} is the dual of homotopy cartesian
(Definition~\ref{defhomotopycartesian}): A square is homotopy
cocartesian if the corresponding square in the opposite category is
homotopy cartesian.

\begin{prop}\label{emprop}
If the square on the left is a homotopy cartesian square of connected
finite type simplicial sets with $K$ simply connected,
\[
\diagram
W\rto\dto&L\dto&&C^{*}W&C^{*}L\lto\\
X\rto&K&&C^{*}X\uto&C^{*}K\uto\lto
\enddiagram
\]
then the square on the right is a homotopy cocartesian square of \einf
algebras. 
\end{prop}

This proposition allows us to prove the following reduction.

\begin{prop}\label{kreduc}
Theorems~\ref{rationalpieceone}, \ref{rationalpiecetwo},
and~\ref{aqvanish} hold for an arbitrary connected finite type
nilpotent simplicial set if they hold for $K(\Z,n)$ and $K(\ZP,n)$ for
all $p$, $n\geq 2$.  
\end{prop}

\begin{proof}
Choose a principally refined Postnikov tower $\{Y_{i}\to
K(G_{i},n_{i}+1)\}$ with each $G_{i}=\Z$ or $\ZP$ for some $p$.  By
induction, we suppose the theorem holds for $Y_{i-1}$.  We choose
cofibrant approximations and cofibrations, making the following
diagram commute.
\[
\diagram
A_{i-1}\dto_{\sim}&B\lcof\rcof\dto_{\sim}&C\dto_{\sim}\\
C^{*}Y_{i-1}&C^{*}K(G_{i},n_{i}+1)\lto\rto&C^{*}L(G_{i},n_{i}+1)
\enddiagram
\]
If we let $A_{i}$ be the pushout $A_{i-1}\amalg_{B}C$, then the previous
proposition tells us that the map $A_{i}\to C^{*}Y_{i}$ is a
quasi-isomorphism.  The squares
\[ \spreaddiagramcolumns{-1pc}
{\diagram
T(A_{i},\Q)\rfib\dfib&T(C,\Q)\dfib&
T(A_{i},\Qhat)\rfib\dfib&T(C,\Qhat)\dfib&
\AQS(A_{i},k)\rfib\dfib&\AQS(C,k)\dfib\\
T(A_{i-1},\Q)\rfib&T(B,\Q)&
T(A_{i-1},\Qhat)\rfib&T(B,\Qhat)&
\AQS(A_{i-1},k)\rfib&\AQS(B,k)
\enddiagram}
\]
are pullbacks of Kan fibrations of Kan complexes.  Here, in the last
square, we choose the augmentations by choosing an augmentation
$A_{i}\to k$ using Lemma~\ref{maptoc} and the given augmentation  on
our given cofibrant approximation.  We have a commutative 
diagram comparing the fibration square defining $Y_{i}$ with the
square on the left and a homotopy commutative diagram comparing the
fibration square defining $\finite{Y}_{i}$ with the square in the
middle.  Inspection of the long exact sequence of homotopy groups
associated to these fibration squares then gives the result.
\end{proof}

To prove the theorems for $K(G,n)$'s, it is convenient to change
coefficients.  For any commutative ring $R$, we can consider the
category $\Ei_{R}$ of \einf $R$-algebras over the operad $(\oE\otimes
R)$.  We have an extension of scalars functor obtained by tensoring an
\einf algebra with the commutative ring $R$.  Extension of scalars is
the left adjoint of the forgetful functor $\Ei_{R}\to \Ei$: We have a
bijection 
\[
\Ei(A,B) \iso \Ei_{R}(A\otimes R,B),
\]
natural in the \einf algebra $A$ and the \einf $R$-algebra $B$.  The
category $\Ei_{R}$ is a model category with fibrations the surjections
and weak equivalences the quasi-isomorphisms.  In particular, the
forgetful functor $\Ei_{R}\to \Ei$ preserves fibrations and
quasi-isomorphisms, and it follows formally that the extension of
scalars functor preserves cofibrations and quasi-isomorphisms between
cofibrant objects.  We can now give the proofs of
Theorems~\ref{rationalpieceone}, \ref{rationalpiecetwo}, and~\ref{aqvanish}.

\begin{proof}[Proof of Theorems~\ref{rationalpieceone}
and~\ref{rationalpiecetwo}]
By Proposition~\ref{kreduc}, it suffices to consider the case when
$Y=K(G,n)$ for $G=\Z$ or $\ZP$ and $n\geq 2$.  Let
$\finite{Y}=K(\finite{G},n)$.  Choose a cofibrant approximation $A\to
C^{*}Y$.  The map $C^{*}(\finite{Y};\Zhat)\to C^{*}(Y;\Zhat)$ is an
acyclic fibration, and so we can choose a lift of the map $A\to
C^{*}Y\to C^{*}(Y;\Zhat)$ to a map $A\to
C^{*}(\finite{Y};\Zhat)$.  Factor the map $A\to C^{*}(\finite{Y};\Zhat)$ as a
cofibration $A\to \finite{A}$ followed by an acyclic fibration
$\finite{A}\to C^{*}(\finite{Y};\Zhat)$.  Then the composite map
$\finite{A}\to C^{*}(Y;\Zhat)$ is a cofibrant approximation.

In the case when $G=\ZP$, let $E$ be the initial object $\oE(0)$ in
$\Ei$, and let $E\to A$ be the unique map; the map $E\otimes \Q\to
A\otimes \Q$ is then a quasi-isomorphism since the cohomology of
$Y=K(G,n)$ is torsion.  In the case when $G=\Z$, let $E$ be the free
\einf algebra on one generator in degree $n$; we choose a map of \einf
algebras $E\to A$ that sends the generator to any cocycle that
represents in cohomology 
the fundamental class of $H^{*}(K(\Z,n))$.  Again the map $E\otimes
\Q\to A\otimes \Q$ is a quasi-isomorphism; this time because the
rational cohomology of $K(\Z,n)$ is the free graded commutative
$\Q$-algebra on the fundamental class and so is the cohomology of
$E\otimes \Q$.  For $R$ containing $\Q$, the map
\begin{multline*}
T(A;R) = \Ei(A,C^{*}(\Delta[\ssdot];R))\iso 
\Ei_{\Q}(A\otimes \Q,C^{*}(\Delta[\ssdot];R))\to\\
\Ei_{\Q}(E\otimes \Q,C^{*}(\Delta[\ssdot];R))\iso
\Ei(E,C^{*}(\Delta[\ssdot];R))=T(E;R)
\end{multline*}
is then a homotopy equivalence by the dual form of
Proposition~\ref{preserveswe}.  On the other hand, we can identify
$T(E;R)$ as follows.  

When $G=\ZP$ and $E=\oE(0)$, we have that $T(E;R)$ is a single point.
Thus, $T(A;R)$ is contractible.  It follows that map $\eta\colon
Y\to T(A;\Q)$ and the map $\finite{\eta}\colon \finite{Y}\to
T(A;\Q)$ are rationalizations.

When $G=\Z$ and $E$ is the free \einf algebra on a generator in degree
$n$ with zero differential, the set of $s$-simplices of $T(E;R)$ is
the set of degree $n$ cocycles of $C^{*}(\Delta[s];R)$; 
in other words, $T(E;R)$ is the standard Eilenberg--Mac\,Lane complex
$K(R,n)$ (see for example \cite[\S 23]{simpob}).  We can assume
without loss of generality that the models for $K(\Z,n)$ and
$K(\Zhat,n)$ we have chosen for $Y$ and $\finite{Y}$ are strictly
$(n-1)$-connected (have only one vertex and no non-degenerate
$i$-simplices for $1\leq i\leq n-1$); for example, the standard models
have this property. Then there is a unique cocycle representing each
cohomology class of $H^{n}Y$ and of $H^{n}(\finite{Y};\Zhat)$.
Looking at the explicit description of the unit of the ($C^{*}$, $T$)
adjunction of Section~\ref{secoutline}, we see that the maps $Y\to
T(E;\Q)$ and $\finite{Y}\to T(E;\Qhat)$ are rationalization, and it
follows that the maps $\eta\colon Y\to T(A;\Q)$ and
$\finite{\eta}\colon \finite{Y}\to T(A;\Qhat)$ are rationalization.
\end{proof}

\begin{proof}[Proof of Theorem~\ref{aqvanish}]
By Proposition~\ref{kreduc}, we are reduced to proving the theorem for
$Y=K(\Z,n)$ or $Y=K(\ZP,n)$ for $n\geq 2$.  We start with $C\to C^{*}Y$ a
given cofibrant approximation, and $C\to k$ a given map of \einf algebras.
If $E$ is any cofibrant \einf $k$-algebra and $E\to C\otimes k$ is a
quasi-isomorphism, then the map
\begin{multline*}
\AQS(C;k)=(\Ei/k)(C,k\oplus C^{*}(\Delta[\ssdot];k[-1]))\\
\iso (\Ei_{k}/k)(C\otimes k,k\oplus C^{*}(\Delta[\ssdot];k[-1]))\to 
(\Ei_{k}/k)(E,k\oplus C^{*}(\Delta[\ssdot];k[-1]))
\end{multline*}
is a homotopy equivalence by the dual form of
Proposition~\ref{preserveswe}; here $\Ei_{k}/k$ denotes the category
of \einf $k$-algebras lying over $k$.   We denote the last simplicial
set in the display above as $\AQS(E \bs k;k)$.  In this notation, it
suffices to show that for some such $E$, $\AQS(E\bs k;k)$ is
contractible. 

In the case when $G=\ZP$ and $p$ is different from the characteristic
of $k$, then $H^{*}(K(\ZP,n);k)$ is trivial and we can take $E$ to be
the initial object $\oE(0)\otimes k$; then $\AQS(E\bs k;k)$ consists
of a single point and is therefore contractible.  We now set $p$ to be
the characteristic of $k$.  In the case when $G=\Z$, the map
$C^{*}(K(\ZPhat,n),k)\to C^{*}(K(\Z,n);k)$ is a quasi-isomorphism.  We
can write $K(\ZPhat,n)$ as the limit of a tower
\[
\cdots \to K(\Z/p^{i}\Z,n) \to \cdots K(\ZP,n)
\]
of principal fibrations $K(\Z/p^{i}\Z,n)\to K(\ZP,n+1)$, and the
canonical map $\Colim H^{*}(K(\Z/p^{i}\Z,n);k)\to
H^{*}(K(\ZPhat,n);k)$ is a quasi-isomorphism.  The arguments for
Propositions~\ref{towerprop} and~\ref{kreduc} then reduce the case
$G=\Z$ to the case $G=\ZP$.

We use the work of \cite[\S 6]{einf}, which constructs an explicit
cofibrant approximation of $C^{*}(K(\ZP,n);k)$ for a field $k$ of
characteristic $p$.  Write $\bE$ for the free functor from
differential graded $k$-modules to \einf $k$-algebras; it takes a
differential graded $k$-module $M$ to the \einf $k$-algebra
\[
\bE M = \myop\bigoplus_{j\geq 0} (\oE(j)\otimes
k)\otimes_{k[\Sigma_{j}]} M^{(j)}.
\]
Let $k[n]$ denote the differential graded $k$-module free on one
generator $x$ in degree $n$, and let $Ck[n]$ be the differential graded
$k$-module free on generators in degrees $n-1$ and $n$ with the
differential taking the lower generator to the higher generator.
In this notation, \cite[\S 6]{einf} describes a map
$\wp\colon \bE k[n]\to \bE k[n]$ such that the pushout over the
inclusion $\bE k[n]\to \bE Ck[n]$,
\[
\diagram
\bE k[n]\dto_{\wp}\rcof&\bE Ck[n]\dto\\
\bE k[n]\rcof&E
\enddiagram
\]
$E=\bE k[n]\amalg_{\wp}\bE Ck[n]$ (which was denoted as $B_{n}$ in
\cite[6.2]{einf}) is quasi-isomorphic to $C^{*}(K(\ZP,n);k)$.  All we
need about the map $\wp$ is that it takes the generator $x$ of $\bE
k[n]$ to the generator $x$ of $\bE k[n]$ minus a class $\rho$ in
\[
(\oE(p)\otimes k)\otimes_{k[\Sigma_{p}]} (k[n])^{(p)}\subset \bE k[n].
\]
(The class $\rho$ represents $P^{0}x$ in $H^{*}(\bE k[n])$.)

Choosing a quasi-isomorphism $E\to C\otimes k$, we obtain an
augmentation $E\to k$.  From the pushout square above, we get a
pullback square.
\[
\diagram
\AQS(E \bs k;k)\rto\dfib&\AQS(\bE Ck[n] \bs k;k)\dfib\\
\AQS(\bE k[n] \bs k;k)\rto_{\wp^{*}}&\AQS(\bE k[n] \bs k;k)
\enddiagram
\]
The vertical arrows are fibrations and the simplicial set in the upper
right hand corner is contractible.  Since $k\oplus
C^{*}(\Delta[\ssdot];k[-1])$ has the square zero multiplication and
the degree of $x-\rho$ is bigger than $0$, any map $\alpha \colon \bE
k[n]\to k\oplus C^{*}(\Delta[\ssdot];k[-1])$ sends $x-\rho$ to the
same element it sends $x$; it follows that the bottom horizontal map
$\wp^{*}$ is an isomorphism.  Thus, $\AQS(E \bs k;k)$ is contractible.
\end{proof}

%%%%%%%%%%%%%%%%%%%%%%%%%%%%%%%%%%%%%%%%%%%%%%%%%%%%%%%%%%%%%%%%%%%%%%%%
\section{Function Complexes and Continuous Functors}\label{secdk2}

The purpose of this section is to set up the machinery we need to
prove Theorems~\ref{looppiece} and~\ref{thma0}.  The main difficulty
is that the identification of $\Tder(C^{*}(Y);\Zhat)$ as 
$\Lambda \finite{Y}$ required comparing with $\Tder(C^{*}(Y);R)$ for
$R=\FP$ and this makes it 
difficult to compare with $\Tder(C^{*}(Y);R)$ for $R=\Q$ and $R=\Qhat$.
The basic idea is to produce a version of $\Lambda
\finite{Y}$ that is a representable functor and to produce a version
of $\Tder(C^{*}Y;R)$ that is a continuous functor.  We then can
able to apply the Yoneda Lemma to construct \emph{natural} maps and to
identify natural transformations.  As a side benefit this gives us
sufficient naturality to prove Theorem~\ref{thmai} in
Section~\ref{secthmai}. 

To carry out this strategy, we use the theory of ``function
complexes'' developed in the papers \cite{dk1,dk2,dk3} of W. G. Dwyer
and D. M. Kan.  This theory works best when we apply it to a small
category.  The set-theoretic technicalities involved in using a
category that is not small are usually treated by ignoring them.  This
is harmless for most applications; essentially the only time these
technicalities become an issue is in the context of mapping space
adjunctions.  Unfortunately, this is the context in which we are
working, and so we are forced to deal with them.

To address these issues, we arrange to work in a small category by
limiting the size of the \einf algebras and simplicial sets we
consider.  For unrelated reasons we explain below, in addition to the
category $\Ei$ of \einf $\Z$-algebras, we also need to work in
the category $\Ei_{\Zhat}$ of \einf $\Zhat$-algebras (over the operad
$\Ei\otimes \Zhat$).  The following definition describes the precise
categories we use.

\begin{defn}
For a simplicial set $Y$, write $\num{Y}$ for the cardinality of the
set of non-degenerate simplices of $Y$, and for a differential graded
module $M$, write $\num{M}$ for the product of the cardinalities of
the modules of all degrees.  Let $\card$ be a cardinal at least as big
as $2^{\continuum}$, where $\continuum$ is the cardinality of the
continuum, and at least as big as $\myop\prod \num{\oE(n)}$.  For a
commutative ring $R$ with $\num{R}\leq\continuum$, let
$\Ecard_{R}$ be a skeleton of the full subcategory of \einf $R$-algebras
$A$ satisfying $\num{A}\leq \card$; we write $\Ecard$ for
$\Ecard_{\Z}$.  Let $\Scard$ be a skeleton of the full subcategory of 
simplicial sets $Y$ satisfying $\num{Y}\leq \card$.
\end{defn}

Recall that a skeleton of a category is a full subcategory with one
object in each isomorphism class; the inclusion of a skeleton is an
equivalence of categories.  We can assume without loss of generality
that $\card^{\continuum}=\card$ (by replacing $\card$ with $2^{\card}$
if necessary); then for any
commutative ring $R$ with cardinality at most $\continuum$ and any
simplicial set with at most $\card$ non-degenerate simplices,
$C^{*}(X;R)$ has cardinality at most $\card$ in each degree.  The
normalized cochain functor with coefficients in $R$ therefore defines
a functor $C^{*}(-;R)\colon \Scard\to \Ecard_{R}$.  In addition, we have
the following homotopy theoretic observations; the proofs are
straightforward cardinality arguments given at the end of this section.

\begin{thm}\label{simpcardmodel}
The category $\Scard$ admits the following closed model structures:
\begin{enumerate}
\item Weak equivalences the usual weak equivalences, cofibrations the
injections, and fibrations the Kan fibrations.
\item Weak equivalences the rational equivalences, cofibrations the
injections, and fibrations defined by the right lifting property.
\item Weak equivalences the finite equivalences (maps that induce
isomorphisms on homology with finite coefficients), cofibrations the
injections, and fibrations defined by the right lifting property.
\end{enumerate}
Moreover, every finite type nilpotent simplicial set is weakly equivalent to
an object in $\Scard$
\end{thm}

\begin{thm}\label{einfcardmodel}
Let $R$ be a commutative ring with $\num{R}\leq \continuum$. The
category $\Ecard$ is a closed model category with weak 
equivalences the quasi-isomorphisms, fibrations the surjections, and
cofibrations as described in \cite[2.4]{einf}.
\end{thm}

The version of the Dwyer--Kan theory we use is called the ``hammock
localization'' \cite[\S 2]{dk2}.  It is defined for any small category with
a suitable notion of weak equivalence, but we only need to apply it to
small closed model categories.  In \cite{dk2}, the hammock
localization is defined to be a simplicial category with the same set
of objects as the original category. The
simplicial sets of maps in this category are virtually never Kan
complexes.  To avoid this inconvenience in the next two sections, we
convert this into a topological category.

\begin{defn}\label{defhamloc}
Let $\Modcat$ be a small closed model category, and let
$L^{H}\subdot\Modcat$ be the simplicial category obtained by hammock
localization as defined in \cite[\S 2]{dk2}.  Let $L\Modcat$ be the
topological category obtained from $L^{H}\subdot\Modcat$ by geometric
realization. 
\end{defn}

As mentioned above, the simplicial category $L^{H}\subdot\Modcat$ and
therefore the topological category $L\Modcat$ has the same object set
as $\Modcat$.  The discrete category $\Modcat$ includes
in $L\Modcat$ as a subcategory, and so a map in $\Modcat$ gives us a
map in $L\Modcat$.  We need the following additional properties of
this localization.

\begin{prop}\cite[3.1]{dk2}
The category $\pi_{0}L\Modcat$ is equivalent to the homotopy category
of $\Modcat$.
\end{prop}

\begin{prop}\cite[3.3]{dk2}\label{wehe}
If $f\colon X\to Y$ is a weak equivalence in $\Modcat$, then for any
object $Z$ in $\Modcat$, the maps
\[
f_{*}\colon L\Modcat(Z,X) \to L\Modcat(Z,Y), \qquad
f^{*}\colon L\Modcat(Y,Z) \to L\Modcat(X,Z)
\]
are homotopy equivalences.
\end{prop}

\begin{prop}\cite[3.4]{dk2}\label{lnat}
A functor $F\colon \Modcat_{1}\to \Modcat_{2}$ that preserves weak equivalences
has a canonical extension to a continuous functor $LF\colon
L\Modcat_{1}\to L\Modcat_{2}$.  If $G\colon \Modcat_{2}\to
\Modcat_{3}$ preserves 
weak equivalences, then $L(G\circ F)=LG\circ LF$.
\end{prop}

We warn the reader that $L$ does not preserve natural transformations;
a natural transformation $F\to G$ generally does not extend to a
natural transformation $LF\to LG$.  See for example \cite[3.5]{dk2}
where the diagram commutes only up to homotopy. In the terminology of
category theory, $L$ is a functor but not a $2$-functor.

The theory of resolutions described in Section~\ref{sechocar} gives a
different way of assigning a simplicial set of maps.  The following
proposition, proved in \cite[6.1]{hh} compares these two
constructions. In the statement, we denote the geometric realization
of a simplicial set by $|\cdot|$.

\begin{prop}\label{comparefc}
Let $X\supdot$ be a cosimplicial resolution of $X$ and let $Y$ be a
fibrant object in $\Modcat$.  The inclusions
\[
|\Modcat(X\supdot,Y)|\to |L\Modcat(X\supdot,Y)| \from
L\Modcat(X,Y)
\]
are homotopy equivalences.
\end{prop}

The dual statement for simplicial resolutions and
cofibrant objects also holds since $(L\Modcat)^{\op}$ and
$L(\Modcat^{\op})$ are isomorphic topological categories.
Alternatively, it is readily deducible from
Propositions~\ref{compareresolutions}, \ref{wehe},
and~\ref{comparefc}.  Combining this with
Proposition~\ref{homotopycartesian}, we get the following proposition,
which requires no fibrancy or cofibrancy assumptions.

\begin{prop}\label{hocatwo}
If the square on the left is homotopy cartesian in $\Modcat$, 
\[
\diagram
W\rto\dto&X\dto&
&L\Modcat(A,W)\dto\rto&L\Modcat(A,X)\dto\\
Y\rto&Z&
&L\Modcat(A,Y)\rto&L\Modcat(A,Z)
\enddiagram
\]
then for any $A$ in $\Modcat$, the square on the right is homotopy
cartesian in $\SP$. 
\end{prop}

We now apply the previous observations to the case at hand, the
model structure on $\Ecard$ and the model structures on $\Scard$.
It is convenient to use $\Scard$ to denote the usual model structure
\ref{simpcardmodel}.(i) on $\Scard$; we write $\Scard_{\Q}$ for the
model structure \ref{simpcardmodel}.(ii) with weak equivalences the
rational equivalences, and $\Scard_{f}$ for the model structure
\ref{simpcardmodel}.(iii) with weak equivalences the finite
equivalences.  The following theorem is an immediate consequence of
Proposition~\ref{comparefc}.  For part (iv), let
$S^{1}$ denote the simplicial model for the circle with one vertex
and one non-degenerate $1$-simplex. 

\begin{prop}\label{completion}
We have natural isomorphisms in $\HoSP$
\begin{enumerate}
\item $|Y| \homeq L\Scard(*,Y)$
\item $\rat{|Y|} \homeq L\Scard_{\Q}(*,Y)$
\item $\finite{|Y|} \homeq L\Scard_{f}(*,Y)$
\item $\Lambda \finite{|Y|}\homeq L\Scard_{f}(S^{1},Y)$
\end{enumerate}
\end{prop}

Thus, we have produced a version of $\Lambda \finite{Y}$ that is a
representable functor on $L\Scard_{f}$.  The following proposition in
particular produces a version of $\Tder(C^{*}Y;R)$ that is a
continuous functor of $L\Scard$; it is an
immediate consequence of the dual statement to
Proposition~\ref{comparefc}. 

\begin{prop}\label{newt}
Let $R$ be a commutative ring with $\num{R}\leq \continuum$.
We have a natural isomorphism in the homotopy category $\HoSP$, 
$|\Tder(A,R)|\homeq L\Ecard(A,R)$.
\end{prop}

Unfortunately, when we apply the Yoneda lemma for natural transformations
out of the representable functor $L\Scard_{f}(S^{1},-)$, we need the
target to be a continuous functor of $L\Scard_{f}$, but $LC^{*}$ is
merely a continuous functor of $L\Scard$ since $C^{*}$ does not
preserve general finite equivalences.  The functor $C^{*}(-;\Zhat)$
does preserve finite equivalences, and here is where we need the
category $\Ecard_{\Zhat}$.  We have in mind $R=\Zhat$, $\W$, or
$\Qhat$ in the following proposition. 

\begin{prop}\label{zhatunit}
Let $R$ be a commutative ring containing $\Zhat$ and satisfying $\num{R}\leq\continuum$. Then the map
\[
\alpha_{R}\colon L\Scard_{f}(*,Y)\overto{LC^{*}(-;\Zhat)}
L\Ecard_{\Zhat}(C^{*}(Y;\Zhat),\Zhat)\overto{i_{*}}
L\Ecard_{\Zhat}(C^{*}(Y;\Zhat),R)
\]
is natural in maps of $Y$ in $L\Scard_{f}$.  The map
\[
\beta_{R}\colon L\Ecard_{\Zhat}(C^{*}(Y;\Zhat),R)
\to L\Ecard(C^{*}(Y;\Zhat),R)
\to L\Ecard(C^{*}(Y),R)
\]
is natural in maps of $Y$ in $\Scard$ (but is not natural for maps of
$Y$ in $L\Scard$ or
$L\Scard_{f}$), and is a homotopy equivalence when $Y$ is finite type.
\end{prop}

\begin{proof}
The first statement follows from Proposition~\ref{lnat}.  In the
second display, the first map is natural in $L\Scard$ (again by
Proposition~\ref{lnat}), and the second map is natural in $\Scard$
(but not natural in $L\Scard$ or $L\Scard_{f}$).
When $Y$ is finite type and $A\to C^{*}Y$ is a cofibrant
approximation, then $A\otimes \Zhat\to C^{*}(Y;\Zhat)$ is a cofibrant
approximation in $\Ei_{\Zhat}$ and the composite map
\[
\Ei_{\Zhat}(A\otimes \Zhat,C^{*}(\Delta[\ssdot];R))\to
\Ei(A\otimes \Zhat,C^{*}(\Delta[\ssdot];R))\to
\Ei(A,C^{*}(\Delta[\ssdot];R)) 
\]
is an isomorphism.  Proposition~\ref{comparefc} then implies
that $\beta_{R}$ is a homotopy equivalence for $Y$ 
finite type. 
\end{proof}

To provide a bridge between the statements in Section~\ref{secoutline}
and the propositions above, we offer the following propositions; they
are immediate consequences of Proposition~\ref{comparefc} and its dual.

\begin{prop}\label{newunitmap}
The following diagram in
$\HoSP$ commutes where the top row is (\ref{uniteta}) and the vertical
maps are the equivalences of Propositions~\ref{completion}
and~\ref{newt}. 
\[
\diagram
|Y|\dline_{\homeq}\rto&|\Tder(C^{*}(Y);\Z)|\dline_{\homeq}\rto
&|\Tder(C^{*}(Y);R)|\dline_{\homeq}\\
L\Scard(*,Y)\rto_{LC^{*}}&L\Ecard(C^{*}(Y),\Z)\rto&L\Ecard(C^{*}(Y),R)
\enddiagram
\]
\end{prop}

For the next proposition, let $Y_{f}$ be a fibrant approximation of
$Y$ in $\Scard_{f}$; then $Y_{f}$ is the Bousfield finite completion of
$Y$ and the map $Y\to Y_{f}$ is a model for the finite completion map
$Y\to \finite{Y}$ when $Y$ is finite type nilpotent.  The map $Y\to
Y_{f}$ induces a quasi-isomorphism $C^{*}(Y_{f};\Zhat)\to
C^{*}(Y;\Zhat)$.  As mentioned above, the following proposition is an
immediate consequence of Proposition~\ref{comparefc} and its dual.
Note that the top composite is $\finite{\eta}$ from
Proposition~\ref{rationalpiecetwo} and the bottom composite is
$\beta_{R}\circ \alpha_{R}$ from Proposition~\ref{zhatunit}.

\begin{prop}\label{lastdiagram}
Let $R$ be as in Proposition~\ref{zhatunit}. The following diagram in $\HoSP$
commutes. 
\[ \spreaddiagramcolumns{-1.25pc}
{\diagram
Y_{f}\dline_{\homeq}\rto
&\Tder(C^{*}(Y_{f};\Zhat);R)\dline_{\homeq}
&\Tder(C^{*}(Y;\Zhat);R)\lto_{\homeq}\rto\ddline_{\homeq}
&\Tder(C^{*}Y;R)\ddline^{\homeq}\\
L\Scard_{f}(*,Y_{f})\rto
&L\Ecard(C^{*}(Y_{f};\Zhat),R)\\
L\Scard_{f}(*,Y)\uto^{\homeq}\rrto
&&L\Ecard(C^{*}(Y;\Zhat),R)\rto\ulto^{\homeq}
&L\Ecard(C^{*}Y,R)
\enddiagram}
\]
\end{prop}

Finally, we close this section with the proof of
Theorems~\ref{simpcardmodel} and~\ref{einfcardmodel}.

\begin{proof}[Proof of Theorems~\ref{simpcardmodel}
and~\ref{einfcardmodel}] 
Using basic cardinal arithmetic, it is easy to see that $\Scard$ is
closed under finite limits and colimits and that $\Ecard_{R}$ is closed
under finite limits.  To see that $\Ecard_{R}$ is closed under finite
colimits, note that the free \einf $R$-algebra on a differential graded
$R$-module $M$, $\bE M$, satisfies 
\[
\num{\bE M}\leq 
\myop\prod (\continuum \cdot \num{\oE(n)} \cdot (\num{M})^{n}).
\]
If $F\colon \oD\to \Ei_{R}$ is a diagram of \einf $R$-algebras, the colimit of
\einf $R$-algebras, $\Colim_{\oD}^{\Ei}F$ may be described as a
differential graded module as the coequalizer
\[
\diagram
\bE(\Colim_{\oD} \bE\circ F)\rto<.5ex>\rto<-.5ex>&
\bE(\Colim_{\oD} F)\rto&\Colim_{\oD}^{\Ei}F
\enddiagram
\]
(where ``$\Colim_{\oD}$'' denotes colimit of differential graded
modules).  Thus, when $\oD$ is a finite diagram and $F$ factors
through $\Ecard_{R}$, the colimit satisfies
$\num{(\Colim_{\oD}^{\Ei}F)}\leq \card$, and so
$\Ecard_{R}$ is closed under finite colimits.  

It remains to see that the factorization axioms hold, since $\Ecard_{R}$ and
$\Scard$ inherit the remaining axioms from the model structure on
$\Ecard_{R}$ and the corresponding model structures on $\Scard$.
In the usual model structure on simplicial sets and the Bousfield
local model structures \cite[\S 10--11]{bousloc} the factorizations
are constructed by the small object argument \cite[II.3.3--4]{quil}.
Let $f\colon X\to Y$ be a map, and consider either of the standard
factorizations $X\to Z\to Y$.
By inspection in the case of~(i) and by~\cite[11.1,11.5]{bousloc} in
the case of~(ii) and~(iii), an easy cardinality argument shows that  
\[
\num{Z}\leq \continuum\cdot\num{X}\cdot\num{Y}
\]
It follows that for a map in $\Scard$ the usual factorizations in
$\SSet$ may be performed in $\Scard$, and so $\Scard$ inherits the
listed closed model structures.  The argument in $\Ecard_{R}$ for
Theorem~\ref{einfcardmodel} is entirely similar.
\end{proof}

%%%%%%%%%%%%%%%%%%%%%%%%%%%%%%%%%%%%%%%%%%%%%%%%%%%%%%%%%%%%%%%%%%%%%%%%
\section{The Proof of Theorems~\ref{looppiece} and \ref{thma0}}\label{secyoneda}

This section is devoted to the proof of Theorems~\ref{looppiece}
and~\ref{thma0}.  
%We use Theorem~\ref{rationalpieceone}, which is
%proved in Section~\ref{secpostnikov}, and Theorem~\ref{wittring},
%which is proved in Sections~\ref{secwittring} and~\ref{secpostnikov}.
The argument is to use the 
basic tools of function complexes described in the previous section to
reinterpret the statements in terms of continuous functors and to
apply the Yoneda Lemma.  We begin by stating the following theorem
that combines and refines Theorems~\ref{looppiece} and~\ref{thma0}.

\begin{thm}\label{refinedtheorem}
There is a natural (in $\Scard$) transformation
\[
\lambda\colon L\Scard_{f}(S^{1},Y)\to L\Ecard(C^{*}(Y);\Zhat)
\]
that is a weak equivalence when $Y$ is connected finite type
nilpotent. Moreover, the following diagrams in $\SP$ commute up to
natural homotopy.
\[ \spreaddiagramcolumns{-0.5pc}
{\diagram
L\Scard_{f}(S^{1},Y)\dto\rto^{\lambda}&L\Ecard(C^{*}(Y),\Zhat)\dto
&L\Scard_{f}(S^{1},Y)\dto\rto^{\lambda}&L\Ecard(C^{*}(Y),\Zhat)\dto\\
L\Scard_{f}(*,Y)\rto_{\gamma_{\W}}&L\Ecard(C^{*}(Y),\W)
&L\Scard_{f}(*,Y)\rto_{\gamma_{\Qhat}}&L\Ecard(C^{*}(Y),\Qhat)
\enddiagram}
\]
Here $\gamma_{\W}=\beta_{\W}\circ \alpha_{\W}$ and $\gamma_{\Qhat}=\beta_{\Qhat}\circ \alpha_{\Qhat}$ for the maps $\alpha_{R}$, $\beta_{R}$ of
Proposition~\ref{zhatunit}. 
\end{thm}

It is convenient to use a formulation of Theorem~\ref{wittring} in
terms of $L\Scard_{f}$ and $L\Ecard$.  The precise formulation we need
is the following proposition.  It is an immediate consequence of
Theorem~\ref{wittring}, Proposition~\ref{lastdiagram}, and the Main
Theorem of \cite{einf}.

\begin{prop}\label{newwittring}
Let $Y$ be a connected finite type nilpotent simplicial set in $\Scard$.
The natural map $\beta_{\W}\circ \alpha_{\W}\colon L\Scard_{f}(*,Y)\to
L\Ecard(C^{*}Y,\W)$ is a weak equivalence.
\end{prop}

Using Propositions~\ref{newunitmap} and~\ref{lastdiagram}, after
passing to the homotopy category and dropping naturality,
Theorem~\ref{refinedtheorem} becomes Theorems~\ref{looppiece}
and~\ref{thma0}.  

We now begin the proof of Theorem~\ref{refinedtheorem}.  Recall from
\cite[II\S5]{serrelocalfields} that there is a unique ring
automorphism $\Phi_{p}$ of the Witt vectors $W(\FPbar)$ that induces
the Frobenius on the residue field $\FPbar$.  We have that the
ring of $p$-adic integers $\ZPhat$ is the subring of $W(\FPbar)$ of elements
$a$ satisfying $\Phi_{p}a=a$.  Another property of the Frobenius is
that for every element $a$ of $W(\FPbar)$ there is some element $x$ in
$W(\FPbar)$ that satisfies $\Phi_{p}x-x=a$.  We let $\Phi$ be the
automorphism of 
$\W$ that performs $\Phi_{p}$ on the factor $W(\FPbar)$ for each $p$.  We then
have that $\Zhat$ is the subring of $\W$ of elements $a$ satisfying
$\Phi a=a$, and for every element $a$ of $\W$, there is some element
$x$ in $\W$ that satisfies $\Phi x -x =a$.  It follows that the square 
\[
\diagram
\Zhat\rto\dto&\W\dto^{(\id,\Phi)}\\
\W\rto_{(\id,\id)}&\W\times \W
\enddiagram
\]
is a pullback square and a simple calculation shows that it is
homotopy cartesian, but note that none of the maps in the square are
surjections.  We can factor the bottom map through a fibration very
easily, as 
\[
\W \to C^{*}(I,\W)\to \W\times \W,
\]
where $I=\Delta[1]$ is the standard $1$-simplex.  The maps are the
induced maps on $C^{*}(-;\W)$ of the projection $I\to *$ and the
inclusions of the vertices $\{0,1\}\to I$.  Let $\Zprime$ be the \einf
$\Zhat$-algebra that makes the following square a pullback.
\[
\diagram
\Zprime\rto\dto&\W\dto^{(\id,\Phi)}\\
C^{*}(I;\W)\rto&\W\times \W
\enddiagram
\]
It follows that the induced map $j\colon \Zhat\to \Zprime$ is a quasi-isomorphism.
As a consequence we get the following proposition.

\begin{prop}\label{zprimeprop}
The natural map 
\[
j_{*}\colon L\Ecard_{\Zhat}(C^{*}(Y;\Zhat),\Zhat)\to
L\Ecard_{\Zhat}(C^{*}(Y;\Zhat),\Zprime)
\]
is a homotopy equivalence.
\end{prop}

The endomorphism $\Phi$ of $\W$ restricts to the identity endomorphism
of $\Zhat$, and so replacing $\W$ with $\Zhat$, the analogue of the square defining
$\Zprime$ as a pullback is the square
\[
\diagram
C^{*}(S^{1};\Zhat)\rto\dto&\Zhat\dto^{(\id,\id)}\\
C^{*}(I;\Zhat)\rto&\Zhat\times \Zhat
\enddiagram
\]
describing $C^{*}(S^{1};\Zhat)$ as a pullback.  The inclusion
$\Zhat\to \W$ and the universal property of the pullback induce the
diagonal maps in the following commutative diagram.
\[ \spreaddiagramcolumns{-2pc}%\spreaddiagramrows{-1pc}
{\diagram
C^{*}(S^{1};\Zhat)\rrto\ddto\drto|>\tip
&&\Zhat\drto\ddto\\
&\Zprime\rrto\ddto
&&\W\ddto\\
C^{*}(I;\Zhat)\rrto\drto
&&\Zhat\times \Zhat\drto\\
&C^{*}(I;\W)\rrto
&&\W\times \W
\enddiagram}
\]
We obtain the
following commutative diagram 
by composing $LC^{*}(-;\Zhat)$ with the maps induced by the diagonal
maps above.
\[ \spreaddiagramcolumns{-3pc}%\spreaddiagramrows{-1pc}
{\diagram
L\Scard_{f}(S^{1},Y)\rrto\ddto\drdashed|>\tip^{\theta}
&&L\Scard_{f}(*,Y)\drto\ddto\\
&L\Ecard_{\Zhat}(C^{*}(Y;\Zhat),\Zprime)\rrto\ddto
&&L\Ecard_{\Zhat}(C^{*}(Y;\Zhat),\W)\ddto\\
L\Scard_{f}(I,Y)\rrto\drto
&&L\Scard_{f}(\{0,1\},Y)\drto\\
&L\Ecard_{\Zhat}(C^{*}(Y;\Zhat),C^{*}(I;\W))\rrto
&&L\Ecard_{\Zhat}(C^{*}(Y;\Zhat),\W\times \W)
\enddiagram}
\]
Note that all the diagonal maps in this diagram are natural for maps of $Y$ in
$L\Scard_{f}$.  By Proposition~\ref{hocatwo} and its dual, both squares
in the diagram above are homotopy cartesian in the category of spaces.
Proposition~\ref{newwittring} implies that when $Y$ is connected
finite type nilpotent, the solid diagonal arrows are homotopy
equivalences.  The dashed arrow $\theta$ is then a weak equivalence; we obtain
the following proposition. 

\begin{prop}\label{firstloopprop}
There is a natural transformation 
\[
\theta \colon L\Scard_{f}(S^{1},Y)\to L\Ecard_{\Zhat}(C^{*}(Y;\Zhat),\Zprime)
\]
of functors $L\Scard_{f}\to \SP$ that is a weak equivalence whenever
$Y$ is connected finite type nilpotent.
\end{prop}

Let $f_{0}$ denote the image under $\theta$ of the identity map of
$S^{1}$.  Let $f_{1}$ be a point in
$L\Ecard(C^{*}(S^{1};\Zhat),\Zhat)$ sent to the component of $f_{0}$
under the map $j_{*}$ in Proposition~\ref{zprimeprop}.  Choose a path $F\colon I\to L\Ecard(C^{*}(S^{1};\Zhat),\Zprime)$
connecting $f_{0}$ and $j_{*}f_{1}$.  The Yoneda lemma then gives us a
natural map 
\[
\zeta \colon L\Scard_{f}(S^{1},Y)\to L\Ecard(C^{*}(S^{1};\Zhat),\Zhat)
\]
and a natural homotopy 
\[
\Theta\colon L\Scard_{f}(S^{1},Y)\times I\to 
L\Ecard(C^{*}(S^{1};\Zhat),\Zprime)
\]
between $\theta$ and $j_{*}\circ \zeta$.  Combining this with
Propositions~\ref{zprimeprop} and~\ref{firstloopprop}, we obtain the
following proposition. 

\begin{prop}\label{secondloopprop}
There is a natural transformation 
\[
\zeta \colon L\Scard_{f}(S^{1},Y)\to L\Ecard_{\Zhat}(C^{*}(Y;\Zhat),\Zhat)
\]
of functors $L\Scard_{f}\to \SP$ that is a weak equivalence whenever
$Y$ is connected finite type nilpotent.
\end{prop}

We now prove Theorem~\ref{refinedtheorem}.

\begin{proof}[Proof of Theorem~\ref{refinedtheorem}]
Since $L\Ecard_{\Zhat}(C^{*}(S^{1};\Zhat),\W)$ and
$L\Ecard_{\Zhat}(C^{*}(S^{1};\Zhat),\Qhat)$ are connected, just as
above, the Yoneda Lemma allows us to choose natural (in $L\Scard_{f}$)
homotopies for the following diagrams in $\SP$ for $R=\W$, $\Qhat$.
\[ 
{\diagram
L\Scard_{f}(S^{1},Y)\dto\rto^(.4){\zeta}
&L\Ecard_{\Zhat}(C^{*}(Y;\Zhat),\Zhat)\dto\\
L\Scard_{f}(*,Y)\rto_(.4){\alpha_{R}}
&L\Ecard_{\Zhat}(C^{*}(Y;\Zhat),R)
\enddiagram}
\]

Let $\lambda =\beta_{\Zhat}\circ \zeta$.  Then $\lambda$ is natural in
$\Scard$ and is a weak equivalence for connected finite type nilpotent
$Y$ by Propositions~\ref{zhatunit} and~\ref{secondloopprop}.
We compose the homotopies above with the maps
$\beta_{\W}$ and $\beta_{\Qhat}$ to get the natural (in $\Scard$)
homotopies for the statement.
\end{proof}

%%%%%%%%%%%%%%%%%%%%%%%%%%%%%%%%%%%%%%%%%%%%%%%%%%%%%%%%%%%%%%%%%%%%%%%%
\section{The Proof of Theorem~\ref{thmai}}\label{secthmai}

In this section we deduce Theorem~\ref{thmai}  from
Theorem~\ref{refinedtheorem}.   The argument consists exclusively of
manipulating homotopy pullbacks.

Define $EY$ to be the homotopy pullback of $L\Ecard(C^{*}Y,\Zhat)\to
L\Ecard(C^{*}Y,\Qhat)$ along $L\Ecard(C^{*}Y,\Q)\to
L\Ecard(C^{*}Y,\Qhat)$, 
and let $DY$ be the homotopy pullback of $\lambda$ along $L\Ecard(C^{*}Y,\Q)\to
L\Ecard(C^{*}Y,\Qhat)$.
\[ \spreaddiagramcolumns{-1.0pc}
{\diagram
EY\rdotted|>\tip\ddotted|>\tip&L\Ecard(C^{*}Y,\Zhat)\dto&
&DY\rdotted|>\tip\ddotted|>\tip&L\Scard_{f}(S^{1},Y)\dto^{\lambda}\\
L\Ecard(C^{*}Y,\Q)\rto&L\Ecard(C^{*}Y,\Qhat)&
&L\Ecard(C^{*}Y,\Q)\rto&L\Ecard(C^{*}Y,\Qhat)
\enddiagram}
\]
Then $E$ and $D$ are functors of $Y$ (in $\Scard$) and we have natural
transformations
\[
L\Ecard(C^{*}Y,\Z)\to EY \from DY
\]
that are weak equivalences when $Y$ is connected finite type
nilpotent.  Let $AY$ be the homotopy pullback of $\gamma_{\Qhat}$
along $L\Ecard(C^{*}Y,\Q)\to L\Ecard(C^{*}Y,\Qhat)$.
\[ 
\diagram
AY\rdotted|>\tip\ddotted|>\tip&L\Scard_{f}(*,Y)\dto^{\gamma_{\Qhat}}\\
L\Ecard(C^{*}Y,\Q)\rto&L\Ecard(C^{*}Y,\Qhat)
\enddiagram
\]
Then $A$ is a functor of $Y$ (in $\Scard$) and we have a natural
transformation $D\to A$ induced by the (second) natural homotopy in
Theorem~\ref{refinedtheorem}.

When we restrict to the case when $Y$ is connected finite type
nilpotent, we get a natural transformation in the homotopy category 
\[
\delta \colon |Y|\homeq L\Scard(*,Y) \to L\Ecard(C^{*}Y,\Z)\to EY \homeq DY \to AY.
\]

We can see from Proposition~\ref{completion}.(iii) and
Theorems~\ref{rationalpieceone} and~\ref{rationalpiecetwo} that the
square defining $AY$ is equivalent in the homotopy category to the
arithmetic square for $Y$, and that $AY$ is abstractly weakly equivalent to
$Y$, but for our argument below we need to
see that this particular map is a weak equivalence.

\begin{thm}\label{deltaequiv}
For $Y$ connected finite type nilpotent, the natural map in the
homotopy category $\delta \colon Y\to AY$ constructed above is a weak
equivalence. 
\end{thm}

\begin{proof}
Since $Y$ and $AY$ are both finite type nilpotent spaces, to see that
$\delta$ is a weak equivalence, it suffices to check that 
it is a weak equivalence after finite completion and that it is a
weak equivalence after rationalization.  To do this, it suffices to check
that the composite maps in the homotopy category 
\[
|Y|\to AY \to L\Scard_{f}(*,Y), \qquad |Y|\to AY \to L\Ecard(C^{*}Y,\Q)
\]
are a finite equivalence and a rational equivalence respectively.  We
can use the fact that the maps comparing $L\Ecard(C^{*}Y,\Z)$, $EY$,
$DY$, and $AY$ are induced by maps of squares to analyze these
composite maps.
The map from $DY\to AY$ sends $L\Scard(S^{1},Y)$ to $L\Scard(*,Y)$ by
the map induced by the inclusion of the vertex in $S^{1}$.  We
therefore see that
that the composite map $|Y|\to L\Scard_{f}(*,Y)$ is a finite equivalence
because the first diagram in Theorem~\ref{refinedtheorem} commutes in
the homotopy category and the composite map $|Y|\to
L\Ecard(C^{*}Y,\W)$ is finite completion by
Proposition~\ref{newwittring}.  Theorem~\ref{rationalpieceone} and
Proposition~\ref{newunitmap} identify the composite map $|Y|\to
L\Ecard(C^{*}Y,\Q)$ as a rational equivalence.
\end{proof}

Theorem~\ref{deltaequiv} is the main result we need for the proof of
Theorem~\ref{thmai}  

\begin{proof}[Proof of Theorem~\ref{thmai}]
For $Y$ connected finite type nilpotent, we define $\epsilon$ to be
the composite of the natural map in the homotopy category 
\[
L\Ecard(C^{*}Y,\Z)\to EY \homeq DY \to AY
\]
and $\delta^{-1}$.  When $Y$ is not connected, \cite[3.1]{finite} allows
us to break up $L\Ecard(C^{*}Y,\Z)$ naturally into a disjoint union of
$L\Ecard(C^{*}Y_{0},\Z)$ over the components $Y_{0}$ of $Y$.  We
then define $\epsilon$ componentwise.
\end{proof}

%%%%%%%%%%%%%%%%%%%%%%%%%%%%%%%%%%%%%%%%%%%%%%%%%%%%%%%%%%%%%%%%%%%%%%%%

%%%%%%%%%%%%%%%%%%%%%%%%%%%%%%%%%%%%%%%%%%%%%%%%%%%%%%%%%%%%%%%%%%%%%%%%
% Bibliography
%%%%%%%%%%%%%%%%%%%%%%%%%%%%%%%%%%%%%%%%%%%%%%%%%%%%%%%%%%%%%%%%%%%%%%%%

\def\MR#1{}
\def\bysame{\leavevmode\hbox to3em{\hrulefill}\thinspace}

\end{document}